\documentclass{conm-p-l}

\usepackage{amssymb}

\def\Ubar{\bar{U}}

\def\hbar{\bar{h}}

\def\iso{\buildrel \sim\over\to}

\def\GS{{\mathfrak{S}}}

\def\Gg{{\mathfrak{g}}}
\def\Gh{{\mathfrak{h}}}
\def\Gm{{\mathfrak{m}}}
\def\Gn{{\mathfrak{n}}}

\def\Gt{{\mathfrak{t}}}

\def\gl{{\mathfrak{gl}}}
\def\Gsl{{\mathfrak{sl}}}
\def\sl{{\mathfrak{sl}}}

\def\CB{{\mathcal{B}}}
\def\CC{{\mathcal{C}}}
\def\CD{{\mathcal{D}}}
\def\CE{{\mathcal{E}}}
\def\CF{{\mathcal{F}}}

\def\CH{{\mathcal{H}}}

\def\CL{{\mathcal{L}}}
\def\CM{{\mathcal{M}}}
\def\CN{{\mathcal{N}}}
\def\CO{{\mathcal{O}}}

\def\CS{{\mathcal{S}}}

\def\CZ{{\mathcal{Z}}}

\def\BA{{\mathbf{A}}}
\def\BB{{\mathbf{B}}}
\def\BC{{\mathbf{C}}}

\def\BF{{\mathbf{F}}}

\def\BH{{\mathbf{H}}}

\def\BP{{\mathbf{P}}}
\def\BQ{{\mathbf{Q}}}
\def\BR{{\mathbf{R}}}

\def\BT{{\mathbf{T}}}

\def\BZ{{\mathbf{Z}}}

\def\Bc{{\mathbf{c}}}

\def\Bk{{\mathbf{k}}}

\def\Bt{{\mathbf{t}}}

\def\eps{\varepsilon}

\def\mcoh{\operatorname{\!-coh}\nolimits}

\def\End{\operatorname{End}\nolimits}

\def\Ext{\operatorname{Ext}\nolimits}

\def\GL{\operatorname{GL}\nolimits}

\def\gr{{\operatorname{gr}\nolimits}}

\def\Hilb{\operatorname{Hilb}\nolimits}

\def\Hom{\operatorname{Hom}\nolimits}

\def\Id{\operatorname{Id}\nolimits}

\def\Ind{\operatorname{Ind}\nolimits}

\def\Irr{\operatorname{Irr}\nolimits}

\def\Lie{\operatorname{Lie}\nolimits}
\def\mMod{\operatorname{\!-mod}\nolimits}

\def\opp{{\operatorname{opp}\nolimits}}

\def\Proj{\operatorname{Proj}\nolimits}

\def\rank{\operatorname{rank}\nolimits}

\def\SL{\operatorname{SL}\nolimits}

\def\Spec{\operatorname{Spec}\nolimits}

\def\Supp{\operatorname{Supp}\nolimits}
\def\Sym{\operatorname{Sym}\nolimits}

\def\Tr{\operatorname{Tr}\nolimits}

\def\ie{{\em i.e.}}

\def\tPhi{{\tilde{\Phi}}}

\newtheorem{thm}{Theorem}[section]

\newtheorem{cor}[thm]{Corollary}
\newtheorem{prop}[thm]{Proposition}
\newtheorem{defi}[thm]{Definition}

\theoremstyle{definition}
\newtheorem{rem}[thm]{Remark}

\newtheorem{ana}{Analogy}
\newtheorem{problem}{Problem}

\usepackage{xy}
\xyoption{all}

\def\KZ{\operatorname{KZ}\nolimits}

\def\eu{\operatorname{eu}\nolimits}
\def\ogr{\operatorname{ogr}\nolimits}
\def\ord{\operatorname{ord}\nolimits}
\def\Sym{\operatorname{Sym}\nolimits}
\def\mqmod{\operatorname{\!-qmod}\nolimits}
\def\mfilt{\operatorname{\!-filt}\nolimits}

\title{Representations of rational Cherednik algebras}
\author{Rapha\"el Rouquier}
\address{Rapha\"el Rouquier\\
Institut de Math\'ematiques de Jussieu --- CNRS\\
UFR de Math\'ematiques, Universit\'e Denis Diderot\\
2, place Jussieu\\
75005 Paris, FRANCE}
\email{rouquier@math.jussieu.fr}

\date{April 2005}
\begin{document}
\begin{abstract}
This paper surveys the representation theory of rational Cherednik algebras.
We also discuss the representations of the spherical subalgebras. We
describe in particular the results on category $\CO$. For type $A$, we
explain relations with the Hilbert scheme of points on $\BC^2$.
We insist on the analogy with the representation theory of complex
semi-simple Lie algebras.
\end{abstract}

\maketitle
{\small
\tableofcontents
}

\section{Introduction}

Let $G$ be a complex reductive algebraic group,
$T$ a maximal torus and $W=N_G(T)/T$ the Weyl group.
Let $\Gt=\Lie T$.

There are several ``Hecke'' algebras associated to $G$ (or to $W$):

$$
\xy
(85,10);(85,-10) **\dir{-} *\dir{>},
(80,20) *{Degeneration},
(-60,-20);(60,-20) **\dir{-} *\dir{>},
(0,-25) *{Affinization},
(0,0) *{\text{
\begin{tabular}{|c|c|c|}
\hline
finite Hecke algebra & affine Hecke algebra &
double affine Hecke algebra \\
 $\CH$ & $\BC[\BT]\otimes\CH$ &
$\BC[\BT]\otimes \BC[\BT^\vee]\otimes\CH$ \\
\hline
& degenerate affine Hecke &
degenerate (trigonometric) daha \\
& $\BC[\Gt]\otimes\BC[W]$ &
 $\BC[\BT]\otimes \BC[\Gt^*]\otimes\BC[W]$ \\
\hline
&&doubly degenerate (rational) daha \\
&& $\BC[\Gt]\otimes \BC[\Gt^*]\otimes\BC[W]$ \\
\hline
\end{tabular}}}
\endxy
$$
\medskip

Here, ``daha'' stands for double affine Hecke algebra (also called Cherednik
algebra) and the structure is given as a $\BC$-vector space. These
algebras have various incarnations:

\medskip
$\bullet\ $The finite Hecke algebra is a quotient of the group algebra of
the braid group, which is the fundamental group of $\Bt_{reg}/W$. There
is a similar description of the affine Hecke algebra (use the space
$\BT_{reg}/W$) as well as of the double affine Hecke algebra \cite{Ch5}.

\smallskip

$\bullet\ $The finite Hecke algebra appears as a coset algebra for a group
of the same
type as $G$, over a finite field. There are similar realizations
of the affine Hecke algebra (use a $1$-dimensional local field) and of the
double affine Hecke algebra ($2$-dimensional local  field, cf \cite{Ka}).

\smallskip
$\bullet\ $There is a geometric realization of (a quotient of) the
daha as the equivariant $K$-theory of the loop Steinberg variety \cite{Va}
(cf also \cite{GaGr,GiKaVa}), generalizing the realization of the affine
Hecke algebra as the equivariant $K$-theory of the ordinary Steinberg
variety.
As pointed out in \cite[\S 7.1]{BeEtGi2}, it is likely that there is an
analogous description of the degenerate daha obtained by using homology
instead of $K$-theory (generalizing the realization of the degenerate
affine Hecke algebra as the equivariant homology of the Steinberg variety).
There is no hint of existence of a geometric realization of the rational daha.

\smallskip
$\bullet\ $When $W$ has type $A_{n-1}$, there are Schur-Weyl dualities between
any of the six types of Hecke algebras and corresponding Lie algebras:
$\Gsl_n$ in the finite case, quantum $\Gsl_n$ and Yangian $\Gsl_n$
in the affine and degenerate affine case,
toroidal quantum $\Gsl_n$ \cite{VarVas1} in the daha case,
toroidal Yangian $\Gsl_n$ in the degenerate daha case, and a
subalgebra of the latter in the rational case \cite{Gu2}.

\medskip
After suitable completions, the daha and the degenerate daha can
be viewed as {\em trivial} deformations of the rational daha
(\cite[p. 283]{EtGi}, \cite[p.65]{Ch}, \cite[\S 7.1]{BeEtGi2}).
As a consequence, the categories of finite dimensional modules for
ordinary, degenerate and rational daha's are equivalent
(\cite{Ch}, \cite[Proposition 7.1]{BeEtGi2}, \cite{VarVas}).
The categories $\CO$ for the
ordinary and degenerate daha's are related \cite{VarVas}. 
Finally, category $\CO$ for the rational daha can be realized as a full
subcategory of category $\CO$ for the degenerate daha \cite{Su}.

\bigskip
Double affine Hecke algebras are related to combinatorics and they were
introduced by Cherednik as a crucial instrument in the proof of the
Macdonald's constant term conjectures. The degenerations (trigonometric
and rational) were obtained in a straightforward way.

In this survey, we are concerned with the representation theory of
rational daha's. Rational Cherednik algebras are connected to
\begin{itemize}
\item
Finite Hecke algebras
\item
Resolutions and deformations of symplectic singularities
\item
Hilbert schemes of points on surfaces.
\end{itemize}
There are also the following connections, which we won't discuss in this
survey:
\begin{itemize}
\item
Integrable systems: rings of quasi-invariants and quantum Calogero-Moser
systems (cf \cite{EtSt} for a survey)
\item
Analytic representation theory, Bessel functions,
unitary representations (cf the book
\cite{Chl}).
\end{itemize}

Many of these aspects actually make sense in the framework of
symplectic reflection algebras \cite{EtGi}.

\medskip
The rational Cherednik algebra is a deformation of the algebra 
$\BC[\Gt\times\Gt^*]\rtimes W$ depending on parameters $t,c$.

The main idea in the study of representations of
rational Cherednik algebras (at $t=1$)
is to handle them like universal enveloping
algebras of semi-simple complex Lie algebras and study in particular a
``category $\CO$''.
We emphasize this in these
notes, by expounding the analogies (cf also the table in \S \ref{table}).
The (finite) Hecke algebra controls a large part of the representation
theory, via the construction of Knizhnik-Zamolodchikov connections.
Another important feature of category $\CO$ is that it generalizes, for
any Weyl group (or even any complex reflection group), the construction of
$q$-Schur algebras.

The usual interplay between representation theory and geometry is incomplete
in general. In type $A$, there are more geometric objects at hand
and the Hilbert scheme of points in $\BC^2$ plays the role of the cotangent
bundle of the flag variety.

The representation theory is quite different for $t=0$. It is related
to (generalized) Calogero-Moser spaces.

\medskip
Section \S\ref{secDunkl} is independent from the rest of the text.
It explains how the rational Cherednik algebras occur naturally in the
study of a commuting family of operators deforming the partial derivatives.

\smallskip
Although the theory has been developed quite intensively over the last
few years, many problems remain and we have listed a number of them.

I have tried to give detailed references for most results, I apologize in
advance for possible omissions.

I thank Ivan Cherednik, Pavel Etingof,
Victor Ginzburg, and Iain Gordon for many useful comments and discussions.
I wish also to thank the mathematics department of Yale University, and in
particular Igor Frenkel, for the invitation to spend the spring semester,
where this paper was written.

\section{A motivation via Dunkl operators}
\label{secDunkl}

\subsection{Dimension $1$}
\label{onedim}
\subsubsection{}
Fix $k\in\BR$. Given $f:\BR\to\BR$ a function of class $C^1$, consider the
function
$$T(f):x\mapsto f'(x)+k\frac{f(x)-f(-x)}{x}.$$
The operator $T$ deforms the ordinary derivation, and presents new
features for special values of $k$.

For example, one shows easily that there exists a non-constant polynomial
killed by $T$ if and only if $k\in -\frac{1}{2}+\BZ_{\le 0}$.

\subsubsection{}
Let us now study the spectrum of $T$. 
We consider the Banach algebra of functions
$B=\{f=\sum_{n\ge 0}a_nX^n:\ ]-1,1[\ \to\BR,\ \sum_n |a_n|<\infty\}$.

We want to solve the equation
\begin{equation}
\label{spectrum}
T(f)=\lambda f \textrm{ and }f(0)=1
\end{equation}
for some $\lambda\in\BR$ and $f\in B$.

Assume $k>0$.
Define
$$\chi:B\to B,\ \ 
f\mapsto (x\mapsto \alpha\int_{-1}^1 f(xt)(1-t)^{k-1}(1+t)^kdt)$$
where $\alpha=\left(\int_{-1}^1 (1-t)^{k-1}(1+t)^kdt\right)^{-1}$
(so that $\chi(1)=1$).
Then, one shows that 
$$T\circ\chi=\chi\circ \frac{d}{dx}.$$
So,
$\chi(\exp(\lambda x))$ is the unique solution of (\ref{spectrum})
(it can be expressed in terms of Bessel functions).

More classical would be the study of the eigenfunctions of the
operator $T^2$ acting on even functions: given $f$ with $f(-x)=f(x)$, then
$T^2(f)=\frac{d^2f}{dx^2}+\frac{2k}{x}\frac{df}{dx}$.

\smallskip
We refer to \cite{ChMa} and \cite[\S 2]{Chl}
for a more detailed study, in particular of the
analytic aspects (Hankel transform, truncated Bessel functions).
\subsection{Dimension $n$}
\subsubsection{}
We are now going to generalize the previous construction to
the case of $n\ge 2$ variables. We will focus on the algebraic aspects
(polynomial functions) and work with complex coefficients. In particular,
we take now $k\in\BC$.

Let $V=\bigoplus_n \BC \xi_i$ and
$V^*=\bigoplus_n \BC x_i$ with the dual basis.
Let $\GS_n$ be the symmetric group on $\{1,\ldots,n\}$. It acts on
$V$ by permutation of the coordinates, hence it acts
on functions $V\to\BC$. We denote by $\rho_{ij}$ the endomorphism
of the space of functions $V\to\BC$ given by the transposition $(ij)$.

Given $1\le i\le n$ and $f:V\to\BC$ smooth, we define
$$T_i(f)=\frac{\partial f}{\partial\xi_i}+k\sum_{j\not=i}
\frac{f-\rho_{ij}(f)}{X_i-X_j}.$$

We have a family of operators (the Dunkl operators)
deforming the ordinary partial derivatives.
What makes this deformation interesting is Dunkl's result:
$$T_i\circ T_j=T_j\circ T_i\textrm{ for all }i,j.$$
Note also that $T_i$ sends a polynomial to a polynomial.

Let $\CE$ be the set of values of $k$ for which there are
non-constant polynomials killed by $T_1,\ldots,T_n$.

The case $n=2$ here is related to the $1$-dimensional case of
\S \ref{onedim} by restricting functions $\BC^2\to\BC$ to the subspace
$x_1+x_2=0$.

Note that the study of functions, in a space similar to $B$ above,
that are simultaneous eigenvectors for all $T_i$'s can be done similarly,
when $k\in\BR_{>0}$. The endomorphism $\chi$ can be constructed first
for polynomials, and then extended to $B$ by continuity. Nevertheless,
there is no explicit integral form for $\chi$.

\subsubsection{}
We denote by $\BC[V]=\BC[X_1,\ldots,X_n]$ the algebra of polynomial
functions on $V$.
Let $H$ be the subalgebra of $\End_{\BC}(\BC[V])$ generated by
$T_1,\ldots,T_n$, $X_1,\ldots,X_n$ (acting by multiplication), and
$\GS_n$ (acting by permutation on the $X_i$'s). This is the
rational Cherednik algebra.

One shows that $k\in\CE$ if and only if $\BC[V]$ is not an irreducible
representation of $H$. Let us now say more about the structure of $H$.

When $k=0$, then $H=D(V)\rtimes\GS_n$, where $D(V)$ is the algebra of
polynomial differential operators on $V$.
In general, there is a vector space decomposition
$$H=\BC[T_1,\ldots,T_n]\otimes \BC[\GS_n]\otimes \BC[X_1,\ldots,X_n].$$
This shows that $H$ (which depends on the parameter $k$) is a deformation
of $D(V)\rtimes\GS_n$.

This is analogous to the decomposition 
$\gl_n(\BC)=\Gn^+\oplus\Gh\oplus\Gn^-$, where $\Gn^+$ (resp. $\Gn^-$)
are strictly upper (resp. lower) triangular matrices and $\Gh$ diagonal
matrices, or rather analogous to the decomposition of the
enveloping algebra (Poincar\'e-Birkhoff-Witt Theorem)
$$U(\gl_n(\BC))=U(\Gn^+)\otimes U(\Gh)\otimes U(\Gn^-).$$
This analogy is a guide for the study of the representation theory of $H$.

First, one defines a category $\CO$ of finitely generated
$H$-modules on which the $T_i$'s act locally nilpotently.
Given $E$ a complex irreducible representation of $\GS_n$, one
gets $\Delta(E)=\Ind_{\BC[T_1,\ldots,T_n]\rtimes \GS_n}^H E$, an object of
$\CO$. It has a unique simple quotient $L(E)$, and one obtains this
way all simple objects of $\CO$. Note that $\Delta(\BC)=\BC[V]$ is the
original faithful representation.
Outside a countable set of values of $k$, then $\CO$ is semi-simple.

\subsubsection{}
We now relate $\CO$ to the Hecke algebra $\CH$
of $\GS_n$, at $q=\exp(2i\pi k)$.

Let $V_{reg}=\{(z_1,\ldots,z_n)\in V| z_i\not=z_j\textrm{ for }i\not=j\}$.
Note that $T_i\in D(V_{reg})\rtimes \GS_n$ and one gets an embedding
$H\subset D(V_{reg})\rtimes \GS_n$. This induces an isomorphism of
algebras
$$H\otimes_{\BC[V]}\BC[V_{reg}]\iso D(V_{reg})\rtimes \GS_n$$
After localization, the deformation becomes trivial!

Let $M\in\CO$. Then, $M\otimes_{\BC[V]}\BC[V_{reg}]$ is
an $\GS_n$-equivariant vector bundle on $V_{reg}$ with a flat connection, that
is shown to have regular singularities (along the hyperplanes and at
infinity). This provides us with a system of differential equations,
and taking solutions we obtain an $\GS_n$-equivariant local system on
$V_{reg}$, \ie, a local system on $V_{reg}/\GS_n$. This corresponds
to a finite dimensional representation of the braid group
$B_n=\pi_1(V_{reg}/\GS_n,(1,2,\ldots,n))$. That representation is shown
to come from a representation of $\CH$ and this defines a functor
$$\KZ:\CO\to\CH\mMod$$
(actually, a contravariant functor; one needs to dualize
or equivalently use the de Rham functor instead of the solution functor
in order to have a covariant functor). This functor
has good homological properties and there is an equivalence
$\CO\simeq \End_\CH(P)\mMod$, where $P$ is a certain $\CH$-module.

When $k\not\in\frac{1}{2}+\BZ$, $P$ can be identified as the
$q$-tensor space $L^{\otimes_q n}$, where $L$ is an $n$-dimensional vector
space, and this identifies $\CO$ with the category of modules
over the $q$-Schur algebra, \ie, a full subcategory of the
category of modules over the quantum general linear group
$U_q(\gl_n)$. The knowledge of character formulas for simple modules in
that setting allows to deduce the multiplicities $[\Delta(E):L(F)]$
for the algebra $H$. This is nicely described in terms of canonical basis
for the Fock space, under the action of $\hat{\sl}_d$, where
$d$ is the order of $k$ in $\BC/\BZ$.

\medskip
Let us finally describe the set $\CE$:
$$\CE=\{-\frac{r}{s}|2\le s\le n,\  r\in\BZ_{>0}\textrm{ and }
(s,r)=1\}.$$
The space of polynomials killed by the $T_i$'s and its structure as
a representation of $\GS_n$ can be determined (cf \cite{DuDeJOp,Du}).

\smallskip
The analytic aspects, which we are not considering here,
are very interesting:
multi-dimensional Bessel functions, Hankel transform
\cite{Op,DuDeJOp}, and
unitary representations of rational Cherenik algebras. It sheds new
light on a number of classical results (cf \cite{Duh}, \cite[Chapter 2]{Chl}).

\section{Structure}
Let us start with the definition of the rational Cherednik algebras and
some important subalgebras, and their main properties,
following \cite[Part 1]{EtGi}.
\subsection{The rational Cherednik algebra}
\label{defrational}
\subsubsection{}
\label{secdef}
Let $W$ be a finite reflection group on a finite dimensional real vector
space $V_\BR$ and let $V=\BC\otimes_\BR V_\BR$. Let $n=\dim_\BC V$.
Let $\CS$ be the set
of reflections of $W$ and $\bar{\CS}=\CS/W$.
Given $s\in \CS$, let $v_s\in V$ (resp.
$\alpha_s\in V^*$) be a $-1$ eigenvector for $s$ acting on $V$ (resp. $V^*$).

Let $\Bc=\{\Bc_s\}_{s\in\bar{\CS}}$ be a family of variables and
$\BA=\BC[\{\Bc_s\},\Bt]$.
Note that for types ADE, we have $|\bar{\CS}|=1$.

\medskip
The {\em rational Cherednik
algebra}
 $\BH$ associated to $(W,V)$ is the
quotient of $\BA\otimes_\BC (T(V\oplus V^*)\rtimes W)$ by the
relations\footnote{In \cite{GGOR}, we use ${\mathbf{\gamma}}_H=-2\Bc_s s$ and
$\Bk_{H,1}=-\Bc_s$, where $H$ is the reflecting hyperplane of $s$}
$$[\xi,\eta]=0 \text{ for } \xi,\eta\in V,\ \
[x,y]=0 \text{ for } x,y\in V^*$$
$$[\xi,x]=\Bt\langle \xi,x\rangle-2\sum_{s\in\CS}
\frac{\langle\xi,\alpha_s\rangle\langle v_s,x\rangle}
     {\langle v_s,\alpha_s\rangle} \Bc_s s$$

There is a filtration on $\BH$ given by
$$ F^0\BH=\BA[W],\ 
F^1\BH=(V\oplus V^*)\otimes_\BC \BA[W]\oplus\BA[W], \textrm{ and }
F^i\BH=(F^1\BH)^i\textrm{ for }i\ge 2.$$
Let $\gr \BH=\bigoplus_{i\ge 0}F^i\BH/F^{i-1}\BH$.
The canonical morphism of $\BA$-modules
$(V\oplus V^*)\otimes_\BC \BA[W]\to \gr\BH$ extends to a surjective
morphism of $\BA$-algebras
$\BA\otimes_\BC (S(V)\otimes S(V^*))\rtimes W\to \gr\BH$. The following result asserts
it is actually an isomorphism. This gives
a triangular decomposition of $\BH$ (Cherednik, \cite[Theorem 1.3]{EtGi}):

\begin{thm}
\label{PBW}
We have a canonical isomorphism of $\BA$-modules
$S(V)\otimes_\BC \BA[W]\otimes_\BC S(V^*)\iso\BH$.
In particular, the canonical map of $\BA$-algebras
$$\BA\otimes_\BC (S(V)\otimes S(V^*))\rtimes W\iso \gr\BH$$
is an isomorphism.
\end{thm}

\begin{proof}[About the proof]
Note that it is enough to prove the
isomorphism after applying $-\otimes_\BA\BC$ for any morphism $\BA\to\BC$,
\ie, for specialized parameters.
Then, one can use the faithful representation by Dunkl operators
(cf \S \ref{Dunkl} and \ref{limitDunkl}).
\end{proof}

\medskip
Given $t\in\BC$ and $c=\{c_s\}_{s\in\bar{\CS}}\in\BC^{\bar{\CS}}$, we put
$H_{t,c}=\BH\otimes_\BA\BC$, where the morphism
$\BA\to\BC$ is given by $\Bt\mapsto t$ and $\Bc_s\mapsto c_s$.
Theorem \ref{PBW} shows that $\BH$ is a deformation of $H_{t,c}$.

\begin{ana}
Let $G$ be a semi-simple complex algebraic group, $T$ a maximal torus,
$B$ a Borel subgroup containing $T$, $U^+$ its unipotent radical and
$U^-$ the opposite unipotent subgroup.
Let $\Gg=\Lie G$, $\Gh=\Lie T$,
$\Gn^+=\Lie U^+$ and $\Gn^-=\Lie U^-$.
We have (Poincar\'e-Birkhoff-Witt Theorem)
$$U(\Gg)=U(\Gn^+)\otimes U(\Gh)\otimes U(\Gn^-).$$
We have a filtration of $U(\Gg)$ given by $F^0U(\Gg)=\BC$,
$F^1U(\Gg)=\Gg\oplus\BC$ and $F^iU(\Gg)=(F^1U(\Gg))^i$.
There is a canonical isomorphism
$S(\Gg)\iso \gr U(\Gg)$.
\end{ana}

\begin{rem}
If $W=W_1\times W_2$ and $V=V_1\oplus V_2$ are compatible decompositions,
then there are canonical isomorphisms
$\BH_1\otimes\BH_2\iso\BH$, where $\BH_i$ is the rational Cherednik algebra of
$(W_i,V_i)$.

Let $\CS'$ be a $W$-invariant subset of $\CS$ such that
$c_s=0$ for $s\in\CS-\CS'$. Let $W'$ be the reflection subgroup of $W$
generated by $\CS'$. Then, there is an embedding $\BH'\subset\BH$ and
$\BH$ is a twisted group algebra of $W/W'$ over $\BH'$.
\end{rem}

\subsubsection{Specializations}
For $t\not=0$, we have
$H_{t,c}\iso H_{1,t^{-1}c}$. We put
$H_c=H_{1,c}$.

\smallskip
Consider the case $c=0$:

\begin{itemize}
\item
We have
$H_{0,0}=S(V\oplus V^*)\rtimes W$.
So, 
$H_{0,\Bc}$ is a deformation of $S(V\oplus V^*)\rtimes W$.

\item
We have $H_{1,0}=D(V)\rtimes W$, where $D(V)$ is the Weyl algebra of
$V$ (algebra of algebraic differential operators on $V$). So, $H_{1,\Bc}$
is a deformation of $D(V)\rtimes W$.
\end{itemize}

\begin{ana}
The parameter space $\BC^{\bar{\CS}}$ corresponds to $\Gh^*/W$ and the analog
of $H_{1,c}$ is $\Ubar_\lambda(\Gg)=U(\Gg)/U(\Gg)\Gm_\lambda$,
where $\lambda\in \Gh^*/W$,
and $\Gm_\lambda$ is the maximal ideal of $Z(U(\Gg))$ image of $\lambda$ by
the canonical isomorphism $\Gh^*/W\iso\Spec Z(U(\Gg))$.

Consider the induced filtration on $\Ubar_\lambda(\Gg)$.
Let $\CN$ be the nilpotent cone of $\Gg$. Then, there is a canonical
isomorphism $\BC[\CN]\iso\gr \Ubar_\lambda(\Gg)$.
\end{ana}

\subsubsection{Fourier transform}
Cf \cite[\S 4,5]{EtGi}.

Fix an isomorphism of $\BC[W]$-modules $F:V\iso V^*$.
This extends to an automorphism $F$ of $\BH_\Bc$ given by
$$V\ni \xi\mapsto F(\xi),\ \ V^*\ni x\mapsto -F^{-1}(x)\ \textrm{ and }\ 
W\ni w\mapsto w.$$

More generally, there is an action of $\SL_2(\BA)$ on $\BH$. The action
of $\left(\begin{matrix} a_{11}&a_{12}\\a_{21}&a_{22}\end{matrix}\right)$ is given by
$$V\ni \xi\mapsto a_{22}\xi + a_{21} F(\xi),\ \
 V^*\ni x\mapsto a_{11}x+a_{12}F^{-1}(x)\ \textrm{ and }\ 
W\ni w\mapsto w.$$

\subsubsection{Twist by characters}
\label{twists}
Let $c\in\BC^{\bar{\CS}}$ and $t\in\BC$. Let $\zeta:W\to\{\pm 1\}$ be
a character. There is an isomorphism of $\BC$-algebras
$$H_{t,c}\iso H_{t,c\zeta},\
V\ni\xi\mapsto \xi,\ \ 
V^*\ni x\mapsto x,\ \ 
W\ni w\mapsto \zeta(w)w.$$

\subsubsection{Deformed Euler vector field and canonical grading}
\label{euler}
We consider in the remaining part of \S \ref{defrational} the algebra
$\BH_{\Bc}=\BH_{1,\Bc}=\BH\otimes_{\BC[\Bt]}\BC[\Bt]/(\Bt-1)$.

Let $B$ be a basis of $V$ and $(b^\vee)_{b\in B}$ be the dual basis.
Let $\eu'=\sum_{b\in B}b^\vee b$ be the ``deformed'' Euler vector field,
$z=\sum_{s\in\CS}\Bc_s(s-1)$ and $\eu=\eu'-z$.
Let $h=\frac{1}{2}\sum_{b\in B}(bb^\vee+b^\vee b)$.
Then, $h=\eu+\frac{1}{2}\dim V-\sum_{s\in\CS}\Bc_s$.
An easy computation in $\BH_\Bc$ shows that
$$[\eu,\xi]=-\xi,\ [\eu,x]=x \textrm{ and }[\eu,w]=0\textrm{ for }
\xi\in V,\ x\in V^* \textrm{ and }w\in W.$$

So, the eigenspace decomposition of $\BH_\Bc$ under the action of $[\eu,-]$
puts a grading on $\BH_\Bc$~: $W$ is in degree $0$, $V^*$ in degree $1$
and $V$ in degree $-1$ (note that this defines a grading on $\BH$ as well,
before specializing $\Bt$ to $1$).

Let $M$ be an $\BH_\Bc$-module. We denote by
$M_\alpha$ is the generalized $\alpha$-eigenspace for $\eu$ acting on $M$.
For certain modules, we have a decomposition 
$M=\bigoplus_{\alpha\in\BC}M_\alpha$ and
this gives a {\em canonical $\BC$-grading} on $M$.

\subsubsection{$\sl_2$-triple}
Cf \cite[\S 4]{De}, \cite[\S 3]{BeEtGi3}.

Let $p:V\times V\to\BC$ and $q:V^*\times V^*\to\BC$ be the
$W$-equivariant perfect pairings induced by $F$. If $\CB$ is orthonormal
for $p$, then $p=\sum_b (b^\vee)^2$ and $q=\sum_b b^2$.
An easy computation shows that
$\langle \frac{1}{2}p,h,-\frac{1}{2}q\rangle$ is an $\Gsl_2$-triple in
$\BH_\Bc$.

\subsection{The spherical subalgebra}
\subsubsection{}
Cf \cite[\S 2]{EtGi}.

Let $e=\frac{1}{|W|}\sum_{w\in W}w$, an idempotent of $Z(\BC[W])$.
Let $\BB=e\BH e$ be the {\em spherical subalgebra} 
of $\BH$. Theorem \ref{PBW} gives a canonical isomorphism
$\BA\otimes_\BC \BC[(V^*\times V)/W]\iso \gr\BB$ (for the induced
filtration on $\BB$).

We have canonical isomorphisms
$B_{1,0}\iso D(V)^W$ and $B_{0,0}\iso\BC[(V^*\times V)/W]$.
So, $\BB$ is a deformation of $D(V)^W$ and of $\BC[(V^*\times V)/W]$.

\begin{ana}
The analogy between $\BB$ and $U(\Gg)$ is more accurate than that with
$\BH$: instead of the smooth orbifold $[(V^*\times V)/W]$, one gets the
singular variety $(V^*\times V)/W$ as corresponding to the nilpotent cone.
Also, for a Cherednik algebra of type $A_1$, then
$B_{1,c}$ is isomorphic to an algebra
$\bar{U}_\lambda(\Gsl_2)$, cf \S \ref{A1sl2}.
\end{ana}

\subsubsection{}
The center and the spherical subalgebra have different behaviours, depending
on $t$.

We have $Z(H_{t,c})=\BC$ if $t\not=0$ \cite[Proposition 7.2]{BrGo}.

The algebra $B_{t,c}$ is commutative if and only if $t=0$
\cite[Theorem 1.6]{EtGi}. In particular, we get a structure of Poisson
algebra on $B_{0,c}$.

We have \cite[Theorem 3.1]{EtGi}:
\begin{thm}[``Satake isomorphism'']
We have an isomorphism
$Z(H_{0,c})\iso B_{0,c},\ z\mapsto ze$.
\end{thm}

The {\em Calogero-Moser space} associated to $W$ is $\CC\CM_c=\Spec Z(H_{0,c})$.
It
is a Gorenstein normal Poisson variety \cite[Theorem 1.5 and Lemma 3.5]{EtGi}
and a symplectic variety when smooth \cite[Theorem 1.8]{EtGi}.

There is an inclusion $S(V)^W\otimes S(V^*)^W\subset Z(H_{0,c})$ and
$Z(H_{0,c})$ is a free $(S(V)^W\otimes S(V^*)^W)$-module of rank $|W|$
\cite[Proposition 4.15]{EtGi}. This gives a finite surjective map
$\Upsilon:\CC\CM_c\to V^*/W\times V/W$.

\subsubsection{}

There is a ``double centralizer Theorem'':
\begin{thm}[{\cite[Theorem 1.5]{EtGi}}]
We have canonical isomorphisms
$\BB\iso\End_\BH(\BH e)$ and
$\BH\iso \End_{\BB^\circ}(\BH e)$.
\end{thm}

The bimodule $H_{t,c}e$ induces a Morita equivalence between
$H_{t,c}$ and $B_{t,c}$ (\ie, 
$H_{t,c}e\otimes_{B_{t,c}}-:B_{t,c}\mMod\to H_{t,c}\mMod$ is an equivalence)
if and only if $H_{t,c}=H_{t,c}eH_{t,c}$.
Cf Theorems \ref{semisimpleMorita} and \ref{shiftA}
for cases of Morita equivalence. 
Note that $H_{0,0}$ is not Morita equivalent to $B_{0,0}$ for $W\not=1$
(the first algebra has finite global dimension while the second one
doesn't).

\section{Representation theory at $t\not=0$}
\subsection{Category $\CO$}
Cf \cite[\S 2]{DuOp},
\cite{Gu}, \cite[\S 2,3]{GGOR}, \cite[\S 2]{BeEtGi3}, \cite[Corollary 6.7]{Gi}.
\subsubsection{Decomposition}

Fix a specialization $H=H_c$ (\ie, $t=1$).
Let $\CO'$ be the category of finitely generated $H$-modules that are
locally finite for $S=S(V)$.

Let $\bar{\lambda}\in V^*/W$, \ie, $\bar{\lambda}$ is a morphism of algebras
$S^W\to \BC$.
Let $\CO_{\bar{\lambda}}$ be the subcategory of objects $M$ in $\CO'$ such that
for any $m\in M$ and any $\xi\in S^W$,
then $(\xi-\bar{\lambda}(\xi))^n\cdot m=0$ for $n\gg 0$.

Then,
$$\CO'=\bigoplus_{\bar{\lambda}\in V^*/W} \CO_{\bar{\lambda}}$$

\subsubsection{Principal block}
We focus our study on $\CO=\CO'_0$ (similar descriptions for
arbitrary $\bar{\lambda}$ have been partially worked out). We write
also $\CO_c$ for the category $\CO$ of $H_c$.

Let $\Irr(W)$ be the set of isomorphism classes of irreducible complex
representations of $W$.
Given $E\in\Irr(W)$, let
$N_E=\sum_{s\in\CS}\frac{\Tr(s|E)}{\dim E}c_s=\frac{\Tr(z|E)}{\dim E}+
\sum_{s\in\CS}c_s$.
We define an order on $\Irr(W)$ by
$E<F$ if $N_F-N_E\in\BZ_{>0}$.

\smallskip
Given $E\in\Irr(W)$, let
$$\Delta(E)=\Ind_{S\rtimes W}^H E$$
where $E$ is viewed as an $(S\rtimes W)$-module with $V$ acting as $0$.

Then, we have
\begin{thm}[{\cite[Theorem 2.19]{GGOR}}]
$\CO$ is a highest weight category with standard objects the
$\Delta(E)$'s.
\end{thm}

\begin{proof}[About the proof]
The approach is similar to the one for affine Lie algebras
and makes crucial use of the canonical grading.
\end{proof}

In particular, $\Delta(E)$ has a unique simple quotient $L(E)$ and
$\{L(E)\}_{E\in\Irr(W)}$ is a complete set of representatives of isomorphism
classes of simple objects of $\CO$.

It follows also that if no two distinct elements of $\Irr(W)$ are comparable,
then $\CO$ is semi-simple.
\begin{cor}
For generic values of $c$, then $\CO$ is semi-simple
and $\Delta(E)=L(E)$ for all $E\in\Irr(W)$.
\end{cor}

The costandard object $\nabla(E)$ is the $H$-submodule of
$\Hom_{S(V^*)\rtimes W}(H,E)$ of elements that are locally nilpotent for 
$S$ (where $E$ is viewed as an $(S(V^*)\rtimes W)$-module with $V^*$ acting as $0$).

\medskip
Simple objects don't have self-extensions:
\begin{prop}[{\cite[Proposition 1.12]{BeEtGi2}}]
We have $\Ext^1_{\CO}(L(E),L(E))=0$ for any $E\in\Irr(W)$.
\end{prop}

\subsubsection{Dualities}
Cf \cite[\S 4]{GGOR}.

Let $M\in\CO$. Let $M'$ be the $k$-submodule of $S(V^*)$-locally
nilpotent elements of $\Hom_\BC(M,\BC)$. 
Consider the anti-involution on $H_c$:
$$\psi:H_c\iso H_c^\opp,\ \ \xi\ni V\mapsto -F(\xi),\ \
V^*\ni x \mapsto -F^{-1}(x),\ \ W\ni w\mapsto w^{-1}.$$
Then $M^\vee=\psi_*M'\in\CO$ and this defines a duality
$$\CO\iso\CO^\opp,\ \ M\mapsto M^\vee.$$
We have $L(E)^\vee\simeq L(E)$ and $\Delta(E)^\vee\simeq\nabla(E)$.

\medskip
The anti-involution
$$H_c\iso H_c^\opp,\ \ \xi\ni V\mapsto -\xi,\ \
V^*\ni x \mapsto x,\ \ W\ni w\mapsto w^{-1}$$
provides $H_c$ with a structure of $(H_c\otimes H_c)$-bimodule and
we obtain a duality
$$R\Hom_{H_c}(-,H_c[n]):D^b(H_c\mMod)\iso D^b(H_c\mMod)^\opp.$$
It restricts to a duality
$$D:D^b(\CO)\iso D^b(\CO)^\opp.$$
We have $D(\Delta(E))^\vee\simeq \nabla(E\otimes\det)$ and
$D(P(E))^\vee\simeq T(E\otimes \det)$ (a tilting module). As a
consequence, $\CO$ is equivalent to its Ringel dual.

\subsubsection{Dihedral groups}
\label{dihedral}
The structure of the standard modules for $W=I_2(d)$ a dihedral
group is given in \cite{Chm}.
Let us explain the results, in the simpler case $d=2m+1$ is odd.
We denote by $\tau_l$ the $2$-dimensional irreducible representation
whose first occurrence in $S(V)$ is in degree $l$ (where $1\le l\le m$).
We assume $c>0$ (cf \S \ref{twists} to deduce the case $c<0$).

$\bullet\ $Assume first $c=\frac{r}{d}$ for some $r\in\BZ_{>0}$,
$d\not|r$.
Let $l\in\{1,\ldots,m\}$ such that $r\equiv \pm l\pmod d$.

We have $L(\rho)=\Delta(\rho)=P(\rho)$ if $\rho\not=\BC,\det,\tau_l$
(they form simple blocks of category $\CO$).

We give now the Loewy series of various modules in $\CO$
(socle and radical series coincide):
$$
\Delta(\BC)=P(\BC)=\begin{matrix}L(\BC)\\L(\tau_l)\end{matrix},\ \
\Delta(\det)=L(\det)=T(\det),\ \
\Delta(\tau_l)=\begin{matrix}L(\tau_l)\\L(\det)\end{matrix}$$

$$P(\det)=T(\tau_l)=\begin{matrix}L(\det)\\L(\tau_l)\\L(\det)\end{matrix},\ \ 
P(\tau_l)=T(\BC)=\begin{matrix}L(\tau_l)\\ L(\det)\oplus L(\BC)\\
L(\tau_l)\end{matrix}$$
The only simple finite dimensional module is $L(\BC)$.

\medskip
$\bullet\ $Assume now $c\in\frac{1}{2}+\BZ_{\ge 0}$.

We have $L(\rho)=\Delta(\rho)=P(\rho)$ if $\rho\not=\BC,\det$
(they form simple blocks of category $\CO$).

We have
$$\Delta(\BC)=P(\BC)=\begin{matrix}L(\BC)\\L(\det)\end{matrix},\ \
\Delta(\det)=T(\det)=L(\det),\ \
P(\det)=T(\BC)=\begin{matrix}L(\det)\\L(\BC)\\L(\det)\end{matrix}.
$$

\medskip
$\bullet\ $If $c\in\BZ_{>0}$ or neither $2c$ nor $dc$ are integers,
then $\CO$ is semi-simple.

\smallskip
The structure is more complicated when $d$ is even, for special values
of the parameter. In particular, finite dimensional modules need not
be semi-simple
(this occurs as well for $W$ of type $D_4$ with parameter $c=\frac{1}{2}$
\cite[Example 6.4]{BeEtGi2}).

\subsubsection{}
In view of the analogy with complex semi-simple Lie algebras, 
finite-dimensional representations are particularly interesting.
A particular class of finite dimensional quotients of $\BC[V]=\Delta(\BC)$
has been studied (``perfect representations''): these
are naturally commutative algebras, and they generalize the Verlinde
algebras \cite{Ch}.

\begin{problem}
\label{pbcatO}
\begin{itemize}
\item
Find the multiplicities $[\Delta(E):L(F)]$. They are known
when $W$ is dihedral (cf \S \ref{dihedral}) and when
has type $A_n$ and $c\not\in\frac{1}{2}+\BZ$ (cf Corollary \ref{decnum}).
\item
Describe the category of finite dimensional representations of $H$.
For which values of $c$ does $H$ have non-zero finite dimensional
representations ?
Cf \cite{Ch,De,BeEtGi2,Go2,ChmEt,Va} for studies of finite dimensional
representations.
These questions are solved in type $A_n$, cf Theorem \ref{finiteA}.
\item
Is $\CO$ Koszul ? If so, is it its own Koszul dual (up to a
change of parameters) ?
\end{itemize}
\end{problem}

\begin{ana}
The category $\CO'$ of finitely generated $U(\Gg)$-modules that are
diagonalizable for $U(\Gh)$ and locally finite for $U(\Gn^+)$ splits into
a sum of subcategories corresponding to a fixed central character. 
The finite dimensional representations are semi-simple and the simple
ones correspond to dominant weights.
The principal
block $\CO$ is a highest weight category with standard objects being Verma modules.
The parametrizing set is $W$, with the Bruhat order. The multiplicities
of simple objects in standard objects are given by evaluation at $1$ of the
Kazhdan-Lusztig polynomials for $W$ (Kazhdan-Lusztig conjecture, proven by
Beilinson-Bernstein and Brylinski-Kashiwara).
The principal block $\CO$ is Koszul and it is
equivalent to its Koszul dual (Beilinson-Ginzburg-Soergel).
\end{ana}

\subsection{Dunkl operators and KZ functor}
\subsubsection{Dunkl operators}
\label{Dunkl}
Cf \cite{Chl}, \cite[\S 4]{EtGi}, \cite[\S 2.2]{DuOp}, \cite[\S 5.2]{GGOR}.

Denote by $\BC$ the trivial representation of $W$. Via the canonical
isomorphism $\BC[V]\iso \Delta(\BC)$, we obtain an action of $H$ on
$\BC[V]$. The action of $W$ is the natural action, the action of $\BC[V]$
is given by multiplication, and the action of $\xi\in V$ is given by
the Dunkl operator (for type $A$, this is the same as \S \ref{secDunkl})
$$T_\xi=\partial_\xi+
\sum_{s\in\CS}\frac{\langle \xi,\alpha_s\rangle}{\alpha_s}c_s(s-1).$$
\noindent
This gives a morphism $\rho:H\to D(V_{reg})\rtimes W$, where
$V_{reg}=V-\bigcup_{s\in\CS}\ker\alpha_s$.

\smallskip
A fundamental property is the faithfulness of that representation
(Cherednik and \cite[Proposition 4.5]{EtGi}):
\begin{thm}
The morphism $\rho$ is injective and induces an isomorphism
$H\otimes_{\BC[V]}\BC[V_{reg}]\iso D(V_{reg})\rtimes W$.
\end{thm}

\begin{proof}[About the proof]
One puts a filtration on $H$ with $W$ and $V^*$ in degree $0$,
and $V$ in degree $1$. Then, $\rho$ is compatible with the filtration
on $D(V_{reg})$ given by the order of differential operators
and the associated graded map is injective.
\end{proof}

\smallskip
Note that $\rho(\eu)=\sum_{b\in\CB}b^\vee b\in D(V)\rtimes W$
is the ordinary Euler vector field.

\begin{rem}
Via the canonical isomorphism $D(V_{reg})^W\iso eD(V_{reg})e, f\mapsto ef$,
the restriction of $\rho$ to $B_{1,c}$ gives an injective morphism
$B_{1,c}\to D(V_{reg})^W$.
\end{rem}

\subsubsection{Knizhnik-Zamolodchikov functor}
Cf \cite[\S 5.3-5.4]{GGOR}.

We are going to associate a vector bundle with a connection on
$V_{reg}$ assciated to an object of $\CO$. 
In the case of a standard object $\Delta(E)$, this is
essentially the Knizhnik-Zamolodchikov-Cherednik connection
\cite{Ch2,Ch4,Op} (cf the affine equation in the trigonometric setting
in \cite{Ch4}).
For generic
values of the parameter, every object of $\CO$ is a sum of $\Delta(E)$'s
and this is used in the constructions below via deformation arguments.

\smallskip
Let $M\in\CO$ and $M_{reg}=\rho_*(M\otimes_{\BC[V]}\BC[V_{reg}])$. This
corresponds to a $W$-equivariant vector bundle on $V_{reg}$ with a flat
connection.  It is shown to have regular singularities, by treating
first the case $M=\Delta(E)$.

Applying the
de Rham functor $\Hom_{\CD(V_{reg})}(\CO_{V_{reg}},-)$ gives a $W$-equivariant
locally constant sheaf on $V_{reg}$. This corresponds to a locally
constant sheaf on $V_{reg}/W$, hence to a finite dimensional representation
$F(M)$
of $B=\pi_1(V_{reg}/W)$ (relative to some base point).
Fixing a base point in $(V_\BR)_{reg}$ provides a description of $W$ as
a finite Coxeter group, with set of simple reflections $\CS_0$. It also
provides an identification of $B$ as the corresponding braid group with set
of generators $\{\sigma_s\}_{s\in \CS_0}$.

Let $\CH$ be the Hecke algebra of $W$ with parameters $\{1,-\exp(2i\pi c_s)\}$,
\ie, the quotient of $\BC[B_W]$ by the relations
$(\sigma_s-1)(\sigma_s+\exp(2i\pi c_s))=0$.

Then, the representation $F(M)$ of $B$ factors through a representation
$\KZ(M)$ of $\CH$. This is proven by first computing
the eigenvalues of monodromy when $M$ is
a standard module and for generic values of the parameter.
A deformation argument shows the result in general.

\smallskip
Let $\CO_{tor}$ be the full subcategory of $\CO$ of objects $M$ such that
$M_{reg}=0$.

\noindent
The main properties of $\KZ$ are given in the following Theorem
\cite[Theorem 5.14, Theorem 5.16, and Proposition 5.9]{GGOR}:
\begin{thm}
\label{thmGGOR}
The functor $\KZ$ is exact and it
induces an equivalence $\CO/\CO_{tor}\iso \CH\mMod$.

Given $M,N\in\CO$, the canonical map
$\Hom_{\CO}(M,N)\to \Hom_{\CH}(\KZ(M),\KZ(N))$ is an isomorphism in the
following cases:
\begin{itemize}
\item when $N$ is projective
\item when $c_s\not\in \frac{1}{2}+\BZ$ for all $s\in\CS$ and
$N$ is $\Delta$-filtered.
\end{itemize}
\end{thm}

\begin{proof}[About the proof]
One shows that $\CO\to\CO/\CO_{tor}$ is fully faithful on projective
objects by using the duality $D$. The fully faithfulness of
$\CO/\CO_{tor}\to\CH\mMod$ is a consequence of the Riemann-Hilbert
correspondence (Deligne). The essential surjectivity follows from a deformation
argument.
The statement in the case where $N$ is $\Delta$-filtered is obtained by
a computation of residues.
\end{proof}

Let $P$ be a progenerator for $\CO$. Then, 
$\CO\simeq \End_{\CH}(\KZ(P))\mMod$. The algebra 
$\End_{\CH}(\KZ(P))$ should be viewed as a ``generalized $q$-Schur algebra''
associated to $W$, cf Theorem \ref{qschur}.

\begin{cor}
The category $\CO$ is semi-simple if and only if $\CH$ is semi-simple.
\end{cor}

\begin{problem}
\begin{itemize}
\item
Provide an explicit construction of a progenerator.
\item
What is the image of a progenerator $P$ of $\CO$ ? What is
$\End_{\CH}(\KZ(P))$ ? This is understood when $W$ has type $A_n$ and
$c\not\in\frac{1}{2}+\BZ$, cf \S \ref{secqschur}.
\end{itemize}
\end{problem}

\begin{ana}
Let $W$ be the Weyl group of $G$ and let
$C=\BC[\Gh]/(\BC[\Gh]\BC[\Gh]_+^W)$ be the algebra of coinvariants.
Note that there is a canonical isomorphism $C\iso H^*(G/B)$.  
Let $P=P(w_0)$ be the ``antidominant'' projective of $\CO$.
There is an isomorphism $C\iso \End_{\CO}(P)$,
the functor $\Hom_{\CO}(P,-):\CO\to C\mMod$ is fully faithful when restricted to
projectives, and the image of a suitable
progenerator is
$$\bigoplus_{w\in W}\BC[\Gh]\otimes_{\BC[\Gh]^{s_1}}
\BC[\Gh]\otimes_{\BC[\Gh]^{s_2}}\cdots\otimes_{\BC[\Gh]^{s_{r-1}}}
\BC[\Gh]\otimes_{\BC[\Gh]^{s_r}}C,$$
where $w=s_1\cdots s_r$ is a reduced decomposition (Soergel).
\end{ana}

\begin{rem}
When $s\mapsto c_s$ is constant (equal parameter case), the $\CH$-modules
$\KZ(\Delta(E))$ are the ``standard'' modules occurring in Kazhdan-Lusztig
theory \cite[Theorem 6.8]{GGOR} (cf \S \ref{secqschur} for type $A$).
\end{rem}

\subsection{Primitive ideals and supports}
Cf \cite[\S 6]{Gi}.
\subsubsection{}
Given $M$ a finitely generated $H$-module, there is a structure of
filtered $H$-module on $M$ such that $\gr M$ is a finitely generated
$\gr H$-module (a ``good filtration''). It is in particular a finitely
generated
$S(V\oplus V^*)$-module (cf \S \ref{secdef}).
Let $\mathrm{Supp}(M)$ be the support of that
module, a $W$-stable closed subvariety of $V^*\times V$.
It is independent of the choice of a good filtration of $M$.

\smallskip
If $M\in\CO$, then $\Supp(M)\subseteq \{0\}\times V$.
When $c=0$, Bernstein's inequality asserts that $\dim\Supp(M)\ge\dim V$.
But there are values of $c$ and objects $M\in\CO_c$
with $\dim\Supp(M)<\dim V$, cf \S \ref{modulesA1}.

We have of course $\Supp(\Delta(E))=\{0\}\times V$, since
the restriction of $\Delta(E)$ to $P$ is free.

\subsubsection{}
Recall that an ideal of $H$ is {\em primitive} if it is the annihilator
of a simple $H$-module. 

\begin{thm}[{\cite[Corollary 6.6]{Gi}}]
\label{duflo}
Every primitive ideal of $H$ is the annihilator of a simple object of $\CO$.
\end{thm}

Let $I$ be an ideal of $H$. Give $I$ the filtration induced by the canonical
filtration on $H$. Then, $\gr I$ is an ideal of
$\gr H=\BC[V^*\times V]\rtimes W$, hence defines a $W$-invariant
closed subvariety of $V^*\times V$, the
{\em associated variety} of $I$.

A {\em parabolic subgroup} of $W$ is the pointwise stabilizer in $W$
of a subspace of $V$. We denote by $\mathrm{Par}(W)$ the set of parabolic
subgroups of $W$.

\begin{thm}[{\cite[Proposition 6.4]{Gi}}]
The associated variety of a primitive ideal of $H$ is of the form
$W\cdot (V^*\times V)^{W'}$ for some $W'\in\mathrm{Par}(W)$. In particular,
its image in $(V^*\times V)/W$ is irreducible.
\end{thm}

\subsubsection{}

\begin{thm}[{\cite[Theorem 6.8]{Gi}}]
Given $M$ simple in $\CO$, there is $W'\in \mathrm{Par}(W)$ such that
$\Supp(M)=\{0\}\times (W\cdot V^{W'})$.
\end{thm}

\begin{problem}
\begin{itemize}
\item
Determine $\Supp L(E)$. This generalizes the problem about finite dimensional
$L(E)$'s (Problem \ref{pbcatO}), they correspond to the case $\Supp L(E)=0$.
\item
Study the order on $\Irr(W)$ defined by $E\prec E'$ if
$\Supp L(E)\subset \Supp L(E')$.
\end{itemize}
\end{problem}

\begin{ana}
The associated variety of a primitive ideal of $U(\Gg)$ is the closure of
a nilpotent class, hence it is irreducible (Borho-Brylinski, Joseph,
Kashiwara-Tanisaki).

Every primitive ideal of $U(\Gg)$ is the annihilator of some
simple object of $\CO'$ (Duflo).

The annihilator of $L(w)$ is contained in the annihilator of $L(w')$
if and only if $w$ is smaller than $w'$ for the ``left cell order''
(Joseph, Vogan).
\end{ana}

\subsection{Harish-Chandra bimodules}

\subsubsection{}
\label{HC}
The definitions and results of this section \ref{HC} follow
\cite[\S 3 and \S 8]{BeEtGi3}.

Let $c,c'\in \BC^{\bar{\CS}}$.

\begin{defi}
A $(H_c,H_{c'})$-bimodule is a
{\em Harish-Chandra bimodule} if it is finitely generated and the action of
$a\otimes 1- 1\otimes a$ is locally nilpotent for every
$a\in S(V)^W\cup S(V^*)^W$.
\end{defi}

Such a bimodule is finitely generated as
a left $H_c$-module, as a right $H_{c'}$-module, as a
$(S(V)^W,S(V^*)^W)$-bimodule and as a $(S(V^*)^W,S(V)^W)$-bimodule.

We denote by $\CH\CC_{c,c'}$ the category of Harish-Chandra
$(H_c,H_{c'})$-bimodules.
The inclusion functor $\CH\CC_{c,c'}\to (H_c\otimes H_{c'}^\opp)\mMod$
has a right adjoint $M\mapsto M_{fin}$.

Given $M\in\CH\CC_{c,c'}$ and
$N\in\CH\CC_{c',c''}$, then $M\otimes_{H_{c'}}N\in\CH\CC_{c,c''}$.

\begin{thm}
Assume $c,c'\in \BZ_{\ge 0}^{\bar{\CS}}$. There is a parametrization
$\{V_{c,c'}(E)\}_{E\in\Irr(W)}$
of the set of isomorphism classes of simple objects of $\CH\CC_{c,c'}$
such that given $E_1,E_2\in\Irr(W)$, then
$$\Hom_\BC(\Delta_{c'}(E_1),\Delta_c(E_2))_{fin}\simeq
\bigoplus_{E\in\Irr(W)}\Hom_{\BC[W]}(E\otimes E_1,E_2)\otimes_\BC V_{c,c'}(E).$$

Furthermore, $E\mapsto V_{c,c}(E)$ extends to an equivalence of
monoidal categories
$$\BC[W]\mMod\iso \CH\CC_{c,c}.$$
\end{thm}
\begin{problem}
Describe the structure of the $2$-category with set of objects $\bar{\CS}$,
$1$-arrows the objects of $\CH\CC_{c,c'}$ and $2$-arrows the morphisms
of $\CH\CC_{c,c'}$.
\end{problem}

\subsubsection{}

\begin{thm}[{\cite[Theorem 3.1]{BeEtGi3}}]
\label{semisimpleMorita}
Assume $\CH$ is semi-simple. Then, $H_c$ is a simple algebra and
$H_c e$ gives a Morita equivalence between $H_c$ and $B_c$.
\end{thm}

\begin{proof}[About the proof]
The key point is Theorem \ref{duflo}. It says in particular
that $H_c e$ gives a Morita equivalence if and only if $e$ kills no
simple object of category $\CO$.
\end{proof}

\smallskip
Let $\eps:W\to\{\pm 1\}$ be a one-dimensional representation of $W$.
Let $e_\eps=\frac{1}{|W|}\sum_{w\in W}\eps(w)w$.
Define $1_\eps:\CS\to\BC,\ s\mapsto 
\begin{cases}
1&\textrm{ if }\eps(s)=-1\\
0&\textrm{ otherwise.}
\end{cases}$

\begin{prop}[{\cite[Proposition 4.11]{BeEtGi3}}]
Assume $\CH$ is semi-simple. Then, the algebras $e_\eps H_ce_\eps$ and
$eH_{c-1_\eps}e$ are isomorphic.
\end{prop}

\begin{thm}[{\cite[Theorem 8.1]{BeEtGi3}}]
\label{Moritaregular}
Assume $\CH$ is semi-simple. Let $m\in \BZ^{\bar{\CS}}$. Then, the algebras
$H_c$ and $H_{c-m}$ are Morita equivalent.
\end{thm}

\begin{problem}[{\cite[Conjecture 8.12]{BeEtGi3}}]
Assume
$\CH$ is semi-simple and let $c'\in\BC^{\bar{\CS}}$. If $H_c$ and
$H_{c'}$ are Morita equivalent, show that there is $\zeta:W\to\{\pm 1\}$ a
character such that $c\zeta -c'\in\BZ^{\bar{\CS}}$.
Cf Theorem \ref{classifA} for a partial answer in type $A$.
\end{problem}

\begin{ana}
Two blocks of category $\CO'$ associated to regular weights are equivalent
via a translation functor.
\end{ana}

\section{Representation theory at $t=0$}
\subsection{General representations}
\subsubsection{Limit Dunkl operators}
\label{limitDunkl}
Following \cite[\S 4]{EtGi}, 
the construction of \S \ref{Dunkl} can be done for the algebra $H_{t,c}$,
$t\not=0$, and it then possible to pass to the limit $t=0$. One obtains
an injective algebra morphism
$$H_{0,c}\hookrightarrow \BC[V^*\times V_{reg}]\rtimes W,\ \
x\mapsto x,\ \ \xi\mapsto \xi+\sum_{s\in\CS}c_s
\frac{\langle\xi,\alpha_s\rangle}{\alpha_s}s,\ \ w\mapsto w.$$

\subsubsection{}
Since $H_{0,c}$ is a finitely generated module over its centre
$Z(H_{0,c})$, it follows that all simple $H_{0,c}$-modules are
finite dimensional. Furthermore, the category of finite dimensional
$H_{0,c}$-modules decomposes into a sum of subcategories according to
the central character (a point of $\CC\CM_c$).

The smoothness of the Calogero-Moser space is related to representation
theory of $H_{0,c}$~:
\begin{thm}[{\cite[Theorem 7.8]{BrGo}, \cite[Theorems 1.7 and 3.7, and
Proposition 3.8]{EtGi},
\cite[Lemma 2.8]{GoSm}}]
\label{smooth}
Let $m\in\CC\CM_c$.
The following assertions are equivalent
\begin{itemize}
\item $m$ is a smooth point of $\CC\CM_c$
\item the Poisson bracket of $\CC\CM_c$ is non-degenerate at $m$
\item there is a unique simple $H_{0,c}$-module with central character $m$
\item the simple $H_{0,c}$-modules with central character $m$ have
dimension $\ge |W|$
\item the simple $H_{0,c}$-modules with central character $m$ are
isomorphic to the regular representation of $W$, as $\BC[W]$-modules.
\end{itemize}
\end{thm}

In particular, if $\CC\CM_c$ is smooth, then its points parametrize
the (isomorphism classes of) simple $H_{0,c}$-modules.

\subsubsection{}
There is a stratification of $\CC\CM_c$ by {\em symplectic leaves}
\cite[\S 3]{BrGo}. Given $I$ a Poisson prime ideal of $B_{0,c}$,
the associated symplectic leaf is
the set of points $m\in\CC\CM_c$ such that $I$ is a maximal Poisson ideal
contained in $m$. There are only finitely many symplectic leaves
\cite[Theorem 7.8]{BrGo}.

\smallskip
The representation theory of $H_{0,c}$ doesn't change inside a symplectic
leaf:
\begin{thm}[{\cite[Theorem 4.2]{BrGo}}]
Let $m,m'\in\CC\CM_c$ be two points in the same symplectic leaf.
Then, $H_{0,c}/H_{0,c}m\simeq H_{0,c}/H_{0,c}m'$.
\end{thm}

One has an irreducibility statement for associated varieties of Poisson
ideals:
\begin{thm}[{\cite[Corollary 3]{Ma}}]
The associated variety in $(V^*\times V)/W$ of a Poisson prime
ideal of $B_{0,c}$ is irreducible.
\end{thm}

\subsection{$0$-fiber}
\subsubsection{}
Cf \cite{Go1}.

Let $I$ be the ideal of $H_{0,c}$ generated by $\BC[V]^W_+$ and
$\BC[V^*]^W_+$ and let $\bar{H}_{0,c}=H_{0,c}/I$.
The $\bar{H}_{0,c}$-modules are $H_{0,c}$-modules whose central character 
is in $\Upsilon^{-1}(0)$ and every simple $H_{0,c}$-module with such a
central character is a $\bar{H}_{0,c}$-module.
The blocks of $\bar{H}_{0,c}$ are given by the central character, \ie,
are in bijection with $\Upsilon^{-1}(0)$~:
$$\bar{H}_{0,c}=\bigoplus_{m\in\Upsilon^{-1}(0)}\bar{H}_{0,c}b_m,$$
where $b_m$ is the primitive central idempotent of $\bar{H}_{0,c}$
corresponding to $m$.
Let $Z_m=\Gamma(\Upsilon^*(0)_m)$, where $\Upsilon^*(0)$ is the
scheme theoretic fiber.

We have a vector space decomposition
$\bar{H}_{0,c}=C\otimes \BC[W]\otimes C'$, where
$C=\BC[V^*]/(\BC[V^*]\BC[V^*]^W_+)$ and 
$C'=\BC[V]/(\BC[V]\BC[V]^W_+)$ are the coinvariant algebras.

\subsubsection{}
The {\em Baby-Verma module} associated to $E\in\Irr(W)$ is
$M(E)=\Ind_{C\rtimes W}^{\bar{H}_{0,c}}E$.
It has a unique simple quotient $L(E)$ and
$\{L(E)\}_{E\in\Irr(W)}$ is a complete set of representatives of isomorphism
classes of simple $\bar{H}_{0,c}$-modules.
Define a map $\Irr(W)\to \Upsilon^{-1}(0),\ E\mapsto m_E$, by the
property that $L(E)$ is in the block $\bar{H}_{0,c}b_{m_E}$.

\smallskip
The blocks corresponding to smooth points can be described precisely:
\begin{thm}[{\cite[Corollary 5.8]{Go1}}]
Let $m\in\Upsilon^{-1}(0)$ be a smooth point of $\CC\CM_c$. Then,
$\bar{H}_{0,c}b_m\simeq \mathrm{Mat}_{|W|}(Z_m)$.

If all points in $\Upsilon^{-1}(0)$ are smooth, then
the canonical map $\Irr(W)\to\Upsilon^{-1}(0)$ is bijective and
$\dim Z_{m_E}=(\dim E)^2$.
\end{thm}

\begin{problem}
\begin{itemize}
\item
Is $\bar{H}_{0,c}$ a symmetric algebra ?
\item
Find the (graded) multiplicities $[M(E):L(F)]$.
In particular, what is the block distribution of the $M(E)$'s ?
\end{itemize}
\end{problem}

\begin{rem}
It is known (\cite[Proposition 16.4]{EtGi}, \cite[Proposition 7.3]{Go1})
that $\CC\CM_c$ is singular
for all values of $c$, when $W$ has type different from $A_n$, $B_n$ and
$D_{2n+1}$. It is conjectured that it will always be singular in type
$D_{2n+1}$ ($n\ge 2$).
\end{rem}

Let us give an example where $E\mapsto m_E$ is not bijective
\cite[\S 7.4]{Go1}.
Let $W=G_2$ and fix a generic value of $c$.
Given $E\not=F\in\Irr(W)$, then $m_E=m_F$ if and only
if $E$ and $F$ are the two distinct $2$-dimensional representations.

\section{Type $A$}
\label{TypeA}
\subsection{Structure}
\subsubsection{}
In this section, let $W=\GS_n$ be the symmetric group on $\{1,\ldots,n\}$
in its permutation representation on $V=\BC^n$.
We consider the canonical bases $(\xi_1,\ldots,\xi_n)$ of $V$ and
$(x_1,\ldots,x_n)$ of $V^*$.
Then, $\BH$ is the $\BC[\Bt,\Bc]$-algebra with generators $\GS_n$,
$x_1,\ldots,x_n$ and $\xi_1,\ldots,\xi_n$ and relations

$[\xi_i,\xi_j]=[x_i,x_j]=0$ for all $i,j$

$wx_iw^{-1}=x_{w(i)}$,
$w\xi_i w^{-1}=\xi_{w(i)}$,

$[\xi_i,x_j]=\Bc\cdot(ij)$ if $i\not=j$ and
$[\xi_i,x_i]=\Bt-\Bc\sum_{k\not=i}(ik)$.

\medskip
One can also consider the action of $W$ on the hyperplane
$V'=\ker(x_1+\cdots+x_n)$. The rational Cherednik algebra of $(W,V')$
is canonically isomorphic to the subalgebra $\BH'$ of $\BH$ generated
by $\xi_i-\xi_j$, $x_i-x_j$ and $W$, for $1\le i,j\le n$ and we
have a decomposition
$\BH=\BH'\otimes \BH^1$, where $\BH^1$ is generated by 
$\xi_1+\cdots+\xi_n$ and $x_1+\cdots+x_n$, and is isomorphic to the first Weyl
algebra.

\smallskip
The algebra $H'_c$ is interesting since it carries finite dimensional
non-zero representations for certain values of $c$, while $H_c$ is the
one that relates most directly to Hilbert schemes of points on the plane.
Note that the categories $\CO$ for $H_c$ and $H'_c$ are canonically
equivalent.
\subsubsection{}
\label{CMA}
The variety $\CC\CM_1$ is isomorphic to the ``usual''
Calogero-Moser space (a smooth symplectic variety)
$$\{(M,M')\in \mathrm{Mat}_n(\BC)\times \mathrm{Mat}_n(\BC) |
 \rank([M,M']+\Id)=1\}/\GL_n(\BC),$$
where $\GL_n(\BC)$ acts diagonally by conjugation
\cite[Theorem 11.16]{EtGi}.

At the level of points, this isomorphism is constructed as follows
\cite[Theorem 11.16]{EtGi}:
let $L$ be a simple representation of $H_{0,1}$.
Fix a basis of the $n$-dimensional space $L^{\GS_{n-1}}$.
The actions of $\xi_n$ and $x_n$ on that space give matrices $M$ and $M'$
such that $\rank([M,M']+\Id)=1$.

The morphism $\Upsilon$ sends $(M,M')$ to the pair of roots of the
characteristic polynomials of $M$ and $M'$.

\begin{rem}[Etingof]
Let $w\in W-\{1\}$. Then, there is $i\in\{1,\ldots,n\}$ such that
$w(i)=j\not=i$. We have
$[\xi_i,x_j(ij)w]=[\xi_i,x_j](ij)w=w$,
hence $w\in [H_{0,1},H_{0,1}]$. It follows that the restriction to
$W$ of a representation of $H_{0,1}$ is a multiple of the regular
representation. This proves the smoothness of $\CC\CM_1$, via
Theorem \ref{smooth}.
\end{rem}

\subsection{Category $\CO$}
\subsubsection{}
\label{secqschur}
Let $q=\exp(2i\pi c)$.
The Hecke algebra $\CH$ of $\GS_n$ with
parameters $(1,-q)$ is the $\BC$-algebra with generators
$T_1,\ldots,T_{n-1}$ and relations
$$T_iT_j=T_jT_i \textrm{ if }|i-j|>1,\ 
T_iT_{i+1}T_i=T_{i+1}T_iT_{i+1}\textrm{ and }
(T_i-1)(T_i+q)=0.$$
We use the standard parametrization of $\Irr(W)$ by partitions of $n$.

Let $\lambda=(\lambda_1\ge\cdots\ge\lambda_r>0)$ be a partition of $n$.
Let $\CH(\lambda)$ be the
subalgebra of $\CH$ generated by $T_1,\ldots,T_{\lambda_1-1},
T_{\lambda_1+1},\ldots,T_{\lambda_1+\lambda_2-1},\ldots$. It is
isomorphic to the tensor product of the Hecke algebras of
$\GS_{\lambda_1},\ldots,\GS_{\lambda_r}$.
Given $\lambda$ a partition of $n$, let $d(\lambda)$ be the number
of $r$-uples $(\beta_1,\ldots,\beta_r)$ whose associated multiset is
that of $\lambda$.
Let $M(\lambda)=\Ind_{\CH(\lambda)}^\CH \BC$,
where $\BC$ is the one-dimensional representations of $\CH$ where the
$T_i$'s act as $1$ and let $M=\bigoplus_\lambda M(\lambda)^{d_\lambda}$.

The {\em $q$-Schur algebra of $\GS_n$} is
$S(n)=\End_{\CH}(M)$.
Note that $S(n)\mMod$ is a highest category with parametrizing set the
set of partitions of $n$.

The $q$-Schur algebra occurs also as a quotient of
the quantum general linear group $U_q(\gl_m)$ for $m\ge n$
(via its action on quantum tensor space $(\BC^m)^{\otimes_q n}$) and
when $q$ is a prime power, as a quotient of the
group algebra of the finite group $\GL_n(\BF_q)$.

\smallskip
The category $\CO$ is described as follows, as conjectured in
\cite[Remark 5.17]{GGOR}:
\begin{thm}[{\cite{Rou}}]
\label{qschur}
Assume $c\not\in \frac{1}{2}+\BZ$. Then,
there is an equivalence
$\CO\iso S(n)\mMod$ making the following diagram commutative
$$\xymatrix{
\CO\ar[rr]^\sim\ar[dr]_{\KZ} && S(n)\mMod\ar[dl]^{M\otimes_{S(n)}-} \\
& \CH\mMod
}$$
and sending $\Delta(\lambda)$ to the standard object of $S(n)\mMod$
associated to $\lambda$ if $c\le 0$ and to the transposed partition
of $\lambda$ if $c>0$.
\end{thm}

\begin{proof}[About the proof]
The proof proceeds by deformation: the parameter ring becomes a discrete
valuation ring and at the generic point the categories are semi-simple.
Then, one shows that the image of the $\Delta$-filtered objects under
the Schur and $\KZ$-functors is a full subcategory closed under extensions.
\end{proof}

In particular, there is a progenerator $P$ of $\CO$ such that
$\KZ(P)=M$. Furthermore, the modules $\KZ(\Delta(\lambda))$ are
the $q$-Specht modules.

\medskip
Assume $c\in\BQ$ and let $d$ be the order of $c$ in $\BQ/\BZ$.

Let $\Sym$ be the space of symmetric functions.
Given $\lambda$ a partition of $n$,
let $s_\lambda$ be the corresponding Schur function.

The Fock space $\Sym$ has a natural
action of the affine Lie algebra $\hat{\Gsl}_d$. There is a
lower canonical basis $\{G^-_\lambda\}_{\lambda\textrm{ a partition}}$
of $\Sym$ \cite{LeThi}. By \cite{VarVas}, the multiplicity
$[\Delta_{S(n)}(\lambda):L_{S(n)}(\mu)]$ is the coefficient of
$G^-_\mu$ in a decomposition of $s_\lambda$ in the lower canonical basis
(a generalisation of the Lascoux-Leclerc-Thibon conjecture on Hecke
algebras, proven by Ariki).
So, we deduce the corresponding result for category $\CO$~:

\begin{cor}
\label{decnum}
Assume $c\not\in \frac{1}{2}+\BZ$. Then,
$[\Delta(\lambda):L(\mu)]$ is the coefficient of
$G^-_\mu$ in a decomposition of $s_\lambda$ in the lower canonical basis.
\end{cor}

\begin{rem}
There is a counterpart of Theorem \ref{qschur} in the trigonometric case
\cite{VarVas}, which builds on an explicit computation of monodromy
(this can't be done in the rational case). It might be possible
to deduce Theorem \ref{qschur} from the trigonometric case by using
\cite{Su}.
\end{rem}

\subsubsection{Finite dimensional representations}
They are completely understood
(\cite[Theorem 1.2]{BeEtGi2}, cf also \cite[\S 7.1]{Ch}):

\begin{thm}
\label{finiteA}
The algebra $H'_c$ has non-zero finite dimensional representations if and only
if $c=\pm\frac{r}{n}$ for some $r\in\BZ_{>0}$ with $(r,n)=1$.
When $c$ takes such a value, all finite dimensional representations
are semi-simple and the only irreducible representation
is $L(\BC)$ when $c>0$ and $L(\det)$ when $c<0$.
\end{thm}

\subsection{Shift functors}
\subsubsection{}
We have $\rho(\delta^{-1}e_{\det} H_{c+1} e_{\det} \delta)=
\rho(eH_ce)$ \cite[Proposition 4.1]{BeEtGi2}.
So, left and right multiplication make
$Q_c^{c+1}=\rho(eH_{c+1}e_{\det} \delta)$ into a $(B_{c+1},B_c)$-bimodule.
Let $S_c=Q_c^{c+1}\otimes_{B_c}-:B_c\mMod\to
B_{c+1}\mMod$ (``Heckmann-Opdam shift functor'').

\medskip
The following result generalizes \cite[Proposition 4.3]{BeEtGi2}.

\begin{thm}[{\cite[Theorem 3.3 and Proposition 3.16]{GoSt1}}]
\label{shiftA}
If $c\in\BR_{\ge 0}$ and $c\not\in \frac{1}{2}+\BZ$, then
\begin{itemize}
\item $eH_c\otimes_{H_c}-:H_c\mMod\to B_c\mMod$ is an equivalence
\item $S_c:B_c\mMod\to B_{c+1}\mMod$ is an equivalence
\item $M\mapsto H_{c+1}e_{\det}\delta\otimes_{B_c}eM:
H_c\mMod\to H_{c+1}\mMod$ is an equivalence.
It restricts to an equivalence
$\CO_c\iso \CO_{c+1}$ sending $\Delta_c(E)$ to $\Delta_{c+1}(E)$.
\end{itemize}
\end{thm}

\begin{proof}[About the proof]
The key point is to show that $H_{c+1}e_{\det}H_{c+1}=H_{c+1}$ and
$H_ceH_c=H_c$ for $c\ge 0$. Let us consider the first equality, the
second one has a similar proof.
If the equality fails, then $e_{\det}$ will kill some simple object of
$\CO$ by Theorem \ref{duflo}. One uses the canonical $\BC$-grading on
$\Delta(\lambda)$, $\lambda$ a partition. One shows that the lowest weight
where the representation $\det$ of $W$ appears in $\Delta(\lambda)$
is strictly larger than the lowest weight where it appears in
$\Delta(\mu)$ whenever $\lambda<\mu$ (this is where the assumption
$c\ge 0$ enters). As a consequence, $\det$
occurs in $L(\mu)$, for any $\mu$, hence $e_{\det} L(\mu)\not=0$ and we deduce
the first equality.
Note that the assumption $c\not\in\frac{1}{2}+\BZ$ comes from the use
of Theorem \ref{thmGGOR}.
To check that $\Delta_c(E)$ goes to $\Delta_{c+1}(E)$, one proves
this after localizing to $V_{reg}$ and then show that the equivalence
$\CO_c\iso\CO_{c+1}$ must preserve the highest weight structure.
\end{proof}

\begin{rem}
Note that, by Theorem \ref{semisimpleMorita} and \S \ref{twists}, we deduce
that $H_c$ and $B_c$ are Morita equivalent for every $c\in\BC$ satisfying
the conditions
$$c\not\in\frac{1}{2}+\BZ\textrm{ and }
c\not\in\{-\frac{m}{d}| m,d\in\BZ,\ 2\le d\le n,\ 0<m<d\}.$$
\end{rem}

\subsubsection{}

There is a Morita equivalence classification of the algebras $H_c$,
when $c$ is not algebraic.
\begin{thm}[{\cite[Theorem 2]{BeEtGi1}}]
\label{classifA}
Let $c\not\in\bar{\BQ}$ and $c'\in\BC$. The algebras $H_c$ and $H_{c'}$
are
\begin{itemize}
\item isomorphic if and only if $c'=\pm c$
\item Morita equivalent if and only if $c\pm c'\in\BZ$.
\end{itemize}
\end{thm}

\begin{proof}[About the proof]
The criterion is obtained by computing the traces
$K_0(H_c)\to HH_0(H_c)$.
\end{proof}

\subsection{Hilbert schemes}
\label{secHilb}
The existence of a link between Hilbert schemes of points on $\BC^2$
and rational Cherednik algebras of type $A$ was pointed out in
\cite{EtGi} and \cite[\S 7.2]{BeEtGi2}.
We describe here some of the constructions and results of \cite{GoSt1,GoSt2}.

\subsubsection{Quantization}
Cf \cite[\S 4--6]{GoSt1}.

Let $\Hilb^n\BC^2$ be the Hilbert scheme of $n$ points in $\BC^2$.
Let $X_n$ be the reduced scheme of $\Hilb^n\BC^2\times_{S^n\BC^2}
\BC^{2n}$ (the isospectral Hilbert scheme). Following Haiman,
we have a diagram

\begin{equation}
\label{Haiman}
\xymatrix{
X_n \ar[rr]^-f \ar[d]_-{\textrm{flat}}^-p && \BC^{2n}=V^*\times V \ar[d] \\
\Hilb^n\BC^2 \ar[rr]_-{\textrm{resolution}}^-\tau &&
 S^n\BC^2=(V^*\times V)/W
}
\end{equation}

Denote by $\CZ_n=\tau^{-1}(0)$ the punctual Hilbert scheme.

\smallskip
Fix $c\in\BR_{\ge 0}$, $c\not\in\frac{1}{2}+\BZ$.

Given $i>j\in\BZ_{\ge 0}$, let
$\CB^{jj}=B_{c+j}$ and
$\CB^{ij}=Q_{c+i-1}^{c+i}\otimes_{B_{c+i-1}}
Q_{c+i-2}^{c+i-1}\otimes\cdots\otimes 
Q_{c+j}^{c+j+1}$.
Let $\CB=\bigoplus_{i,j\ge 0}\CB^{ij}$.
This is a non-unital algebra. We denote by $1_i$ the unit of $\CB^{ii}$.
We denote by $\CB\mMod$ the category of finitely generated $\CB$-modules
$M$ such that $M=\bigoplus_{i\ge 0} 1_i M$.
We denote by $\CB\mqmod$ the abelian category quotient of $\CB\mMod$ by
the Serre subcategory of objects $M$ such that $1_iM=0$ for
$i\gg 0$.

The filtration  by the order of differential operators on
$D(V_{reg})\rtimes W$ induces a filtration on $\CB$ and we denote
by $\ogr \CB=\bigoplus_{i\ge j\ge 0}\ogr \CB^{ij}$ the associated graded 
(non-unital) algebra. We define a category $(\ogr \CB)\mqmod$ as above.

\begin{thm}[{\cite[Theorem 6.4]{GoSt1}}]
\label{mainGS}
There are equivalences
$$B_c\mMod\iso \CB\mqmod$$
$$\Hilb^n\BC^2\mcoh\iso (\ogr \CB)\mqmod.$$
\end{thm}

\begin{proof}[About the proof] The first equivalence is an
immediate consequence of the fact that the $\CB^{ij}$'s induce Morita
equivalences (Theorem \ref{shiftA}).

Consider $A^1=\BC[\BC^{2n}]^{\det}$, the $\det$-isotypic part for the action
of $W$ on $\BC[\BC^{2n}]$ and let $A^{\prime 1}$ be the ideal of
$\BC[\BC^{2n}]$ generated by $A^1$.
Let $A^d=(A^1)^d$ and $A^{\prime d}=(A^{\prime 1})^d$ (inside $\BC[\BC^{2n}]$)
for $d\ge 1$ and let $A^0=\BC[\BC^{2n}]^{\GS_n}$ and
$A^{\prime 0}=\BC[\BC^{2n}]$.
Let $A=\bigoplus_{d\ge 0}A^d$ and $A'=\bigoplus_{d\ge 0}A^{\prime d}$.

Then, there are isomorphisms
$X_n\iso \Proj A'$ and $\Hilb^n\BC^2\iso \Proj A$ so that the
diagram (\ref{Haiman})
above becomes the following diagram, with canonical maps (Haiman)
$$\xymatrix{
\Proj \bigoplus_{d\ge 0}A^{\prime d} \ar[rr] \ar[d] &&
  \Spec A^{\prime 0} \ar[d] \\
\Proj \bigoplus_{d\ge 0}A_d \ar[rr] && \Spec A_0
}$$
Now, the Theorem follows from the equalities
$\ogr \CB^{ij}=eA^{i-j}\delta^{i-j}e$ between subspaces of
$\BC[\BC^{2n}]^{\GS_n}$, whose delicate proof involves understanding the
graded structure of these two subspaces.
\end{proof}

\subsubsection{Localization}
Cf \cite{GoSt2}.

Consider the order filtration on $H_c$~: $\ord^0 H_c=\BC[V]\rtimes W$,
$\ord^1 H_c=V\cdot\ord^0 H_c+\ord^0 H_c$ and $\ord^i H_c=(\ord^1 H_c)^i$ for
$i\ge 2$. This induces a filtration on $B_c$.

Theorem \ref{mainGS} gives a functor
$\Phi_c$ from the category $B_c\mfilt$ of $B_c$-modules with a good
filtration (for the order filtration) to the category
$\Hilb^n\BC^2\mcoh$ of coherent sheaves over $\Hilb^n\BC^2$.

We have
$$\Phi_c(B_c)\simeq \CO_{\Hilb^n\BC^2}$$

The image of $eH_c$ with the order filtration is the Procesi bundle
\cite[Theorem 4.5]{GoSt2}:
$$\Phi_c(eH_c)\simeq p_*\CO_{X_n}.$$

We have \cite[Proposition 5.4]{GoSt2} (this relates to 
\cite[Conjectures 7.2 and 7.3]{BeEtGi2}):
$$\Phi_{1/n}(L_{1/n}(\BC))\iso \CO_{\CZ_n}.$$

\medskip
Given $M\in B_c\mfilt$, there is an induced tensor product filtration on
$S_c(M)=Q_c^{c+1}\otimes_{B_c}M$.
Let $\CL$ be the determinant of the
universal rank $n$ vector bundle on $\Hilb^n\BC^2$ (an ample line bundle).
The geometric importance of $S_c$ is
provided by the next result, which explains the constructions of
Gordon-Stafford.

\begin{thm}[{\cite[Lemma 4.4]{GoSt2}}]
There is an isomorphism of functors
$B_c\mfilt\to \Hilb^n\BC^2\mcoh$:
$$\Phi_{c+1}\circ S_c(-)\iso \CL\otimes \Phi_c(-).$$
\end{thm}

\begin{problem}[{\cite[Question 1]{GoSt2}}]
Let $M$ be an $H_c$-module with a good filtration. Then, 
$\gr M$ is a finitely generated $(\BC[V^*\times V]\rtimes W)$-module.
Let $\tPhi(M)\in D^b(\Hilb^n\BC^2\mcoh)$ be its image under the equivalence
of derived categories (Bridgeland-King-Reid, Haiman):
$$(p_*Lf^*(-))^W:D^b((\BC[V^*\times V]\rtimes W)\mMod)\iso
D^b(\Hilb^n\BC^2\mcoh).$$
Is there an isomorphism $\tPhi(M)\iso \Phi(eM)$ ?
Gordon and Stafford construct a morphism which is surjective on $\CH^0$.
A related question is to understand which
$(\BC[V^*\times V]\rtimes W)$-modules
can be quantized to $H_c$-modules for some value of $c\in\BR_{\ge 0}$.
\end{problem}

\begin{rem}
Gordon and Stafford actually work with $H'_c$ and they consider 
the $W$-Hilbert scheme $\Hilb(n)$ of $(V')^*\times V'$. There is an isomorphism
$\Hilb^n\BC^2\iso \Hilb(n)\times \BC^2$, so the geometric properties of
$\Hilb(n)$ and $\Hilb^n\BC^2$ are easily related \cite[\S 4.9]{GoSt1}.
\end{rem}

\subsubsection{Characteristic cycles}
Cf \cite[\S 6]{GoSt2}.

Let $Z=Z(n)=\tau^{-1}(\{0\}\times V/W)$.
Let $\lambda=(\lambda_1\ge\cdots\ge\lambda_r>0)$ be a partition of $n$ and
$S^\lambda\BC^2$ be the subvariety of $S^n\BC^2$
of $0$-cycles of $\BC^2$ of the form
$\sum_i \lambda_i x_i$, where $x_1,\ldots,x_r\in\BC^2$ are distinct.

Let $Z_\lambda$ be the closure of $Z\cap\tau^{-1}(S^\lambda\BC^2)$. This
is a Lagrangian subvariety of $\Hilb^n\BC^2$ and the $Z_\lambda$'s,
where $\lambda$ runs over the partitions of $n$, are the irreducible
components of $Z$ (Grojnowski, Nakajima).

Given $\lambda$ a partition of $n$, let 
$m_\lambda\in\Sym$ be the corresponding monomial symmetric function.
There is an isomorphism (Grojnowski, Nakajima)
$$\xi:\bigoplus_{n\ge 0}H_n(Z(n))\iso \Sym,\ \
[Z_\lambda]\mapsto m_\lambda$$
where $H_n$ is the Borel-Moore homology with complex coefficients.

\medskip
Recall that the cycle support of a coherent sheaf $\CF$ is the cycle
$\sum_i n_i[Z_i]$, where $Z_i$ runs over irreducible components of the support
of $\CF$ and $n_i$ is the dimension of $\CF$ at the generic point of $Z_i$.

Let $M\in\CO_c$. Fix a good filtration of $eM$ and consider the
part of the cycle support of $\Phi(eM)$ involving only subvarieties of
dimension $n$. This is independent of the choice of the good filtration and
this gives an isomorphism \cite[Corollary 6.10]{GoSt2}
$$\gamma:K(\CO_c)\otimes_\BZ\BC\iso H_n(Z)$$

The following Theorem describes the characteristic cycle of $\Delta(\mu)$.

\begin{thm}[{\cite[Theorem 6.7]{GoSt2}}]
Let $\mu$ be a partition of $n$.
The support of $\Phi(e\Delta(\mu))$ is a union of $Z_\lambda$'s.
We have
$$\xi\gamma([\Delta(\mu)])=s_\mu.$$
\end{thm}

\begin{problem}[{\cite[Question 4.9 and \S 6.8]{GoSt2}}]
From Corollary \ref{decnum}, one deduces that $\xi\gamma([L(\mu)])$ is the
lower canonical basis element of the Fock space corresponding to $\mu$. 
Are the irreducible components of
the support of $\Phi(eL(\mu))$ all of dimension $n$ ? If so,
this would completely describe the characteristic cycle of $L(\mu)$.
\end{problem}

\begin{problem}[{\cite[Problem 7.7]{GoSt2}}]
Gordon and Stafford show \cite[Lemma 7.7]{GoSt2} that the top
Borel-Moore homology of $\Hilb^n\BC^2\times_{S^n\BC^2}\Hilb^n\BC^2$ with
the convolution product is isomorphic to the representation ring of
$W$. The determination of the $(\BC^\times)^2$-equivariant $K$-theory
ring under convolution remains to be done.
\end{problem}

\begin{rem}
There is a different geometric approach started in \cite{GaGi}, which also
leads to the construction of characteristic cycles on the Hilbert scheme.
\end{rem}

\subsubsection{}
Let us summarize
$$\xymatrix{
\CB \ar@{.>}[d]_{\textrm{quantization}} \ar@{<.>}[rr]^{\sim} &&
B_c \ar@{.>}[d]^{\textrm{quantization}} \\
\Hilb^n\BC^2 \ar[rr]_{\textrm{resolution}} && (V^*\times V)/W
}$$

The algebra $H_c$ is a ``quantization'' of the orbifold $[(V^*\times V)/W]$.
Furthermore, $X_n$, viewed as a scheme over $\Hilb^n\BC^2$ and
over $[(V^*\times V)/W]$, becomes after ``quantization''
the $(B_c,H_c)$-bimodule $eH_c$.

The transform with kernel $\CO_{X_n}$ gives
an equivalence of triangulated categories (``McKay correspondence'')
$D^b(\Hilb^n\BC^2\mcoh)\iso D^b([(V^*\times V)/W]\mcoh)$.
This is simpler in the non-commutative case, where $eH_c$ gives
an equivalence of abelian categories
$B_c\mMod\iso U_c\mMod$ (for suitable $c$'s).

\begin{ana}
Let $\CB=G/B$ be the flag variety of $G$. Given $w\in W$,
we put $\CB_w=BwB/B$.
We have the Springer resolution of singularities
$T^*\CB\to\CN$.
We have a canonical isomorphism
$\Ubar_0(\Gg)\iso \Gamma(\CB,\CD_\CB)$
and we have
$\Gamma(\CB,\CD_\CB)\mMod\simeq \CD_\CB\mcoh$ (Beilinson-Bernstein).
We have $\gr\CD_\CB\iso\CO_{T^*\CB}$.
$$\xymatrix{
\CD_{\CB} \ar@{.>}[d]_{\textrm{quantization}} \ar@{<.>}[rr]^{\sim} &&
\bar{U}_0(\Gg) \ar@{.>}[d]^{\textrm{quantization}} \\
T^*\CB \ar[rr]_{\textrm{resolution}} && \CN 
}$$

Given an object $M$ of $\CO'$, we can consider the characteristic
cycle of the corresponding $\CD_\CB$-module, an element of
$\bigoplus_{w\in W}\BZ [\overline{T^*_{\CB_w}\CB}]$.

There is a canonical isomorphism of algebras
between the top homology of the Steinberg
variety $Z=T^*\CB\times_\CN T^*\CB$ and $\BC W$ (Kazhdan-Lusztig).
The images of the components
of $Z$ give a basis $(b_w)_{w\in W}$ of $\BC W$. Define
$(\beta_{w,w'})_{w,w'}$ by
$w=\sum_{w'}\beta_{w,w'}b_{w'}$. Then,
the characteristic cycle of $\Delta(w)$ is
$\sum_{w'}\beta_{w,w'}[\overline{T^*_{\CB_{w'}}\CB}]$ (Kashiwara-Tanisaki).

The $(G\times\BC^\times)$-equivariant $K$-theory of $Z$ is isomorphic to the
affine Hecke algebra of type $W$ (Kazhdan-Lusztig, Chriss-Ginzburg).
\end{ana}

\begin{problem}
Can one deform the category of $W$-equivariant mixed Hodge modules on
$V$ ?
\end{problem}

\section{Type $A_1$}
\subsection{Presentation}
Take $W=A_1=\langle s\rangle$ acting on $V=\BC\xi$ and put $V^*=\BC x$ where
$\langle \xi,x\rangle=1$.
Then, $\BH$ is  the $\BC[\Bt,\Bc]$-algebra with generators
$s,x,\xi$ and relations
$$s^2=1,\ sxs=-s,\ s\xi s=-\xi,\ [\xi,x]=\Bt-2\Bc s.$$
We have $e=\frac{1}{2}(1+s)$.

Fix $t,c\in\BC$. Recall that
$H_{t,c}\simeq H_{1,t^{-1}c}$ if $t\not=0$ and
$H_{0,c}\simeq H_{0,1}$ if $c\not=0$. So, there are three types of
algebras in the family~: $H_{1,c}$, $H_{0,1}$ and $H_{0,0}$.

\subsection{Category $\CO$ and $\KZ$}
\label{modulesA1}
\subsubsection{}
We take $t=1$.

We identify $\Delta(\BC)$ with $\BC[x]$. The action of the
generators is given by
\begin{align*}
x: & x^i\mapsto x^{i+1} \\
s: & x^i\mapsto (-1)^i x^i \\
\xi: & x^i\mapsto 
\begin{cases}
i x^{i-1} \textrm{ if }i \textrm{ is even}\\
(i-2c) x^{i-1} \textrm{ if }i \textrm{ is odd}
\end{cases}
\end{align*}
In particular, $\Delta(\BC)=L(\BC)$ if and only if
$c\not\in \frac{1}{2}+\BZ_{\ge 0}$.
If $c=\frac{1}{2}+n$ with $n\ge 0$, then
$x^{2n+1}\BC[x]$ is the radical of $\Delta(\BC)$ and
$L(\BC)=\BC[x]/(x^{2n+1})$.

The Dunkl operator is $T_\xi=\frac{d}{dx}+\frac{c}{x}(s-1)$.
The connection on the trivial line bundle over $V_{reg}$ is 
$\frac{d}{dx}$. The solutions are constant functions, the monodromy
operator is trivial.

\subsubsection{}

We identify $\Delta(\det)$ with $\BC[x]$. The action of the
generators is given by
\begin{align*}
x: & x^i\mapsto x^{i+1} \\
s: & x^i\mapsto (-1)^{i+1} x^i \\
\xi: & x^i\mapsto 
\begin{cases}
i x^{i-1} \textrm{ if }i \textrm{ is even}\\
(i+2c) x^{i-1} \textrm{ if }i \textrm{ is odd}
\end{cases}
\end{align*}
In particular, $\Delta(\det)=L(\det)$ if and only if
$c\not\in -(\frac{1}{2}+\BZ_{\ge 0})$.
If $c=-(\frac{1}{2}+n)$ with $n\ge 0$, then
$x^{2n+1}\BC[x]$ is the radical of $\Delta(\det)$ and
$L(\det)=\BC[x]/(x^{2n+1})$.

The connection on the trivial line bundle over $V$ is
$\frac{d}{dx}+\frac{2c}{x}$. In a neighborhood of $x=1$, we have the
solution $f=x^{-2c}$ with $f(1)=1$. Analytical continuation in a tubular
neighborhood of the path $t\in[0,1]\mapsto \exp(i\pi t)$
gives a function with value $\exp(-2i\pi c)$ at $-1$.
The action of the monodromy operator is $-\exp(-2i\pi c)$.
The $\KZ$-construction involves the de Rham functor while here
we considered the solution functors (=horizontal sections). We pass
from one to the other by dualizing the representation of the 
fundamental group. So, the generator $\sigma_s$ of the braid group
acts now by $-\exp(2i\pi c)$ on $\KZ(\Delta(\det))$.

\subsubsection{}
When $c\not\in\frac{1}{2}+\BZ$, then $\CO$ is semi-simple.

When $c\in -(\frac{1}{2}+\BZ_{\ge 0})$, then we have
an exact sequence $0\to \Delta(\BC)\to\Delta(\det)\to L(\det)\to 0$.
We have $P(\det)=\Delta(\det)$ and
$P(\BC)=H\otimes_{\BC[\xi]\rtimes W}S(\xi)/(\xi)^2$.

When $c\in \frac{1}{2}+\BZ_{\ge 0}$, then we have
an exact sequence $0\to \Delta(\det)\to\Delta(\BC)\to L(\BC)\to 0$.
We have $P(\BC)=\Delta(\BC)$ and
$P(\det)=H\otimes_{\BC[\xi]\rtimes W}(S(\xi)/(\xi)^2\otimes \det)$.

So, when $c\in\frac{1}{2}+\BZ$, then $\CO$ is equivalent to the
principal block of category $\CO$ for $\Gsl_2(\BC)$ (this is no
miracle, cf \S \ref{A1sl2}).

\subsubsection{}
We have $ex+xe=2x$, $e\xi+\xi e=2\xi$, so
$\bigoplus_{i,j}\BC x^i\xi^j\subset HeH$.
Since $[\xi,x]=t-2cs$, we have $t-2cs\in HeH$. So,
$2c+t\in HeH$.
It follows that $HeH=H$ if $2c+t\not=0$.

Assume now $2c=-t$. We put a structure of $H$-module on $L=\BC$ by
letting $x$ and $\xi$ act by $0$ and $s$ by $-1$.
Then, $e$ acts by $0$, so $HeH$ annihilates $L$, hence $HeH\not=H$.
So, we have $HeH=H$ if and only if $2c\not=-t$.
Note that when $c=-\frac{1}{2}$ and $t=1$,
via the identification $\Delta(\det)=\BC[x]$, then
$L\simeq L(\det)=\BC[x]/x\BC[x]$.
Note also that the algebra $B_{1,-\frac{1}{2}}$ is simple.

Finally, the following assertions are equivalent:
\begin{itemize}
\item $H_{t,c}eH_{t,c}=H_{t,c}$
\item $H_{t,c}$ and $B_{t,c}$ are Morita equivalent
\item $c\not=-\frac{1}{2}t$.
\end{itemize}

\subsection{Spherical subalgebra}
\subsubsection{}
One has $ex^i\xi^j e=\begin{cases}
ex^i\xi^j & \textrm{ if }i+j\text{ is even}\\
0 & \textrm{ otherwise}
\end{cases}$
and $ex^i\xi^j e=\pm ex^i s\xi^j e$.
So, $\BB$ has a basis $(ex^i \xi^j)_{i,j\ge 0,i+j\textrm{ even}}$.
It is generated as a $\BC[\Bt,\Bc]$-algebra by $u=\frac{1}{2}e x^2 e$,
$v=-\frac{1}{2}e \xi^2 e$ and $w=e x\xi e$.
We obtain an isomorphism between $\BB$ and the $\BC[\Bt,\Bc]$-algebra
generated by $u$, $v$ and $w$ with the relations
$w(w-\Bt-2\Bc)=-4uv$ and $[u,v]=\Bt w+\Bt(\frac{1}{2}\Bt-\Bc)$.

\subsubsection{}
We have $\BC[u,v,w]/(w(w-2c)+4uv)\iso B_{0,c}$. The variety
$\Spec B_{0,c}$ is smooth if and only if $c\not=0$.

\subsubsection{}
\label{A1sl2}
For $t=1$, we have $[u,v]=w+\frac{1}{2}-c$,
$[v,w]=2v$ and $[u,w]=-2u$.
Let $e_+$, $e_-$ and $h$ be the standard generators of $\Gsl_2(\BC)$ and
let $C=e_+e_-+e_-e_++\frac{1}{2}h^2$ be the Casimir element.
We have an isomorphism \cite[Proposition 8.2]{EtGi}
$$U(\Gsl_2)/\langle C-\frac{1}{2}(c-\frac{1}{2})(c+\frac{3}{2})\rangle\iso
B_{1,c},\ \ e_+\mapsto u,\ \ e_-\mapsto v.$$

We have $B_{1,c}$ Morita equivalent to $B_{1,c+1}$ if and only if
$c\not=-\frac{3}{2},-\frac{1}{2}$.

\subsubsection{}
The representation theory of $U(\Gsl_2)$ is related to
$T^*\BP^1$, since the flag variety of $\Gsl_2$ is a projective line.
Note that $\Hilb^W(V^*\times V)\simeq T^*\BP^1$, the minimal resolution of
singularities of the cone $\Spec B_{0,0}$.

\subsection{Double affine Hecke algebra}
We finish with some words on the daha and its relation with
the rational daha.
\subsubsection{}
The double affine Hecke algebra $\BH^{ell}$ is the
$\BC[\tau^{\pm 1},q^{\pm 1}]$-algebra with generators
$X^{\pm 1},Y^{\pm 1},T$ and relations
\begin{equation}
\label{ellA1}
(T-\tau)(T+\tau^{-1})=0,\ TXT=X^{-1},\ TY^{-1}T=Y \textrm{ and }
Y^{-1}X^{-1}YXT^2=q.
\end{equation}

There is a triangular decomposition
$\BH^{ell}=\BC[X^{\pm 1}]\otimes \CH\otimes \BC[Y^{\pm 1}]$
\cite[\S 1.4.2]{Chl}.
Let us consider the $\BH^{ell}$-module
$\Ind_{\BC[Y^{\pm 1}]\otimes\CH}^{\BH^{ell}}\BC$, where
$Y$ and $T$ act on $\BC$ by multiplication by $\tau$
\cite[Proof of Lemma 1.4.5]{Chl}.
We identify this module with $\BC[X^{\pm 1}]$. This gives a faithful
representation of $\BH^{ell}$. The action of $\BH^{ell}$ is given by
\begin{align*}
X :\ & X^i\mapsto X^{i+1} \\
T :\ & X^i\mapsto X^{-i}T+(\tau^{-1}-\tau)(X^{i-2}+X^{i-4}+\cdots+X^{-i}) \\
YT^{-1} :\ & X^i\mapsto q^{-i}X^{-i}
\end{align*}
So, $T$ acts by $\tau s+\frac{\tau-\tau^{-1}}{X^2-1}(s-1)$ and
$Y$ acts by $sp T$, where $p(f)(X)=f(q^{-1}X)$.

\subsubsection{}
We follow \cite[Proposition 4.10]{EtGi}.

Fix $c\in\BC$.
We consider the ring $\BC[[h]]$ with its $h$-adic topology.
Let $\hat{H}^{ell}$ be the $\BC[[h]]$-algebra topologically
generated by three elements $s$, $x$ and $y$ with relations
(\ref{ellA1}) for $X=e^{hx}$, $Y=e^{hy}$, $T=se^{h^2cs}$, $q=e^{h^2}$ and
$\tau=e^{h^2c}$.

The first relation, taken at order $0$, gives $(s-1)(s+1)=0$.
The second and third relations, taken at order $1$, give
$sxs=-x$ and $sys=-y$. Finally, the last relation, taken at order $2$, gives
$yx-xy=1-2cs$. This gives rise to an isomorphism of $\BC$-algebras
$H\iso \hat{H}^{ell}\otimes_{\BC[[h]]}
\BC[[h]]/(h)$. Note finally that $\hat{H}^{ell}$ is actually a trivial
deformation of $H$, \ie, there is an isomorphism of topological
$\BC[[h]]$-algebras $H\hat{\otimes}\BC[[h]]\iso \hat{H}^{ell}$
\cite[p.65]{Ch}.

\section{Generalizations}
\subsection{Complex and symplectic reflection groups}
\subsubsection{}
The definition of the rational Cherednik algebra generalizes to the case
where $W$ is a complex reflection group on $V$ and more
generally, when $W$ is a symplectic reflection group on a space $L$
(which is $V\oplus V^*$ in the complex reflection case).
Theorem \ref{PBW} remains valid in that setting.
The construction and the main results
on category $\CO$ generalize to complex reflection groups.

Such deformations have been introduced and studied before in
\cite{CrBoHo} in the case of a symplectic reflection group acting
on a symplectic space of dimension $2$.

These constructions of symplectic reflection algebras 
are special cases of a more general construction \cite{Dr}. An even more general
construction is given in \cite{EtGaGi}, where finite groups are
replaced by reductive groups.

Another direction of generalization is globalization, where $V$ is
replaced by an algebraic variety acted on by a finite group \cite{Et}.

\subsubsection{}
Many results of \S \ref{TypeA} should generalize to the
complex reflection groups $B_n(d)\simeq \BZ/d\wr\GS_n$
(cf \cite{Mu} for a generalization of some of the constructions of
\S \ref{secHilb}).
Some geometric aspects (should) generalize even to
$\Gamma\rtimes \GS_n$, where $\Gamma$ is a finite subgroup of $\SL_2(\BC)$.
For example, the Hilbert scheme to consider is
$\Hilb^n (\Hilb^\Gamma\BC^2)$.
The description of multiplicities for $B_n(d)$ should generalize 
Corollary \ref{decnum} using
suitable canonical bases of higher level Fock spaces \cite{Rou}.

Some finite dimensional representations have been constructed in
the case $\Gamma\rtimes \GS_n$. Those where the $x_i$'s and $\xi_i$'s act
by zero have been classified in \cite{Mo1}. More general finite
dimensional simple representations have be shown to exist by cohomological
methods \cite{EtMo,Mo2}. These representations come from ones where
the $x_i$'s and the $\xi_i$'s act by zero via reflection functors
and every simple finite dimensional representation for parameters
``close to $0$'' is of this form \cite{Ga}.

\subsubsection{}
Assume $W=\Gamma\wr\GS_n$ acting on $\BC^{2n}$,
for some finite subgroup $\Gamma$ of $\SL_2(\BC)$,

There is a crepant (=symplectic) resolution $X_c\to \CC\CM_c$
and an equivalence \cite[Theorem 1.2]{GoSm}
$$D^b(H_{0,c}\mMod)\iso D^b(X_c\mcoh)$$
The variety $X_c$ is constructed as a moduli space of representations
of $H_c$ and the kernel of the equivalence is the universal bundle.
The case $c=0$ is the McKay correspondence.

The algebra $H_{0,c}$ is actually a non-commutative crepant resolution of
$\CC\CM_c$, in the sense of Van den Bergh \cite[Lemma 3.10]{GoSm}.

Note that $\CC\CM_c$
is smooth for generic values of $c$ \cite[Proposition 11.11]{EtGi} and
it has then a description generalizing that of
\S \ref{CMA} \cite[Theorem 11.16]{EtGi}.

\subsection{Characteristic $p>0$}
The rational Cherednik algebra can be defined over $\BZ$, and in
particular over an algebraic closure $\bar{\BF}_p$ of the field with
$p$ elements. The representation theory in characteristic $p>0$ is
quite different from that in characteristic $0$, due to the presence
of a big center: we have $\gr Z(H_c)\iso (\BC[V^*\times V]^p)^W$ when $p>n$
(Etingof, cf \cite[Theorem 10.1.1]{BeFiGi}). In particular, 
$H_c$ is a finite dimensional module over its center, hence all its
simple representations are finite-dimensional.

There is a localization result in characteristic $p$.
\begin{thm}[{\cite[Theorem 7.3.2]{BeFiGi}}]
Assume $c\le 0$ and $c\not\in\frac{1}{2}+\BZ$. Assume $p$ is large enough.

Then, there is a sheaf $\CF_c$ of Azumaya algebras over the Frobenius twist
$\Hilb^{(1)}$ of $\Hilb^n\BA^2$ and an equivalence 
$$D^b(\CF_c\mMod)\iso D^b(H_c\mMod).$$
This comes from an isomorphism
$H^0(\Hilb^{(1)},\CF_c)\iso H_c$ and from
the vanishing $H^{>0}(\Hilb^{(1)},\CF_c)=0$.
\end{thm}

Cf \cite{La} for the determination of irreducible representations
of Cherednik algebras associated to the rank $1$ groups $W=\BZ/d$, over a
field of characteristic $p\not|d$, making explicit results of \cite{CrBoHo}.

\section{Table of analogies}
\label{table}
$$
\begin{array}{l|l}
\BH \textrm{ or } \BB & U(\Gg) \\
H_{t,c}=S(V)\otimes\BC[W]\otimes S(V^*) & U(\Gg)=U(\Gn^+)\otimes U(\Gh)
 \otimes U(\Gn^-)\\
\BC^{\bar{\CS}} & \Gh^*/W \\
H_{1,c} \textrm{ or }B_{1,c} & U(\Gg)/U(\Gg)\Gm_{\lambda} \\
(V^*\times V)/W & \CN \\
\textrm{parabolic subgroups of }W & \textrm{nilpotent classes}\\
\Irr(W) & W \\
\CH & S(\Gh^*)/S(\Gh^*)S(\Gh^*)_+^W\simeq H^*(G/B) \\
? & \textrm{dominant weights} \\
\Hilb^n\BC^2 \textrm{ (type $A_{n-1}$)} & T^*G/B \\
\end{array}
$$


\begin{thebibliography}{[GiKapVas]}
\bibitem[BerEtGi1]{BeEtGi1} Y.~Berest, P.~Etingof and V.~Ginzburg,
	{\em Morita equivalence of Cherednik algebras},
	 J. Reine Angew. Math. {\bf 568} (2004), 81--98.
\bibitem[BerEtGi2]{BeEtGi2} Y.~Berest, P.~Etingof and V.~Ginzburg,
	{\em Finite dimensional representations of rational Cherednik
	algebras},
	Int. Math. Res. Not. {\bf 19} (2003), 1053--1088. 
\bibitem[BerEtGi3]{BeEtGi3} Y.~Berest, P.~Etingof and V.~Ginzburg,
	{\em Cherednik algebras and differential operators on
	quasi-invariants},
	Duke Math. J. {\bf 118} (2003), 279--337.
\bibitem[BezFiGi]{BeFiGi} R.~Bezrukavnikov, M.~Finkelberg and V.~Ginzburg,
	{\em Cherednik algebras and Hilbert schemes in characteristic p},
	with appendices by P.~Etingof and V.~Vologodsky,
	preprint math.RT/0312474(v4).
\bibitem[BrGo]{BrGo} K.A.~Brown and I.~Gordon,
	{\em Poisson orders, representation theory and symplectic reflection
	algebras}, J. Reine angew. Math. {\bf 559} (2003), 193--216.
\bibitem[Che1]{Ch1} I.~Cherednik,
	{\em Calculation of the monodromy of some $W$-invariant local
	systems of type $B,C$ and $D$},
	Funct. Anal. Appl. {\bf 24} (1990), 78--79.
\bibitem[Che2]{Ch2} I.~Cherednik,
	{\em Generalized braid groups and local $r$-matrix systems},
	Soviet Math. Dokl. {\bf 40} (1990), 43--48.
\bibitem[Che3]{Ch3} I.~Cherednik,
	{\em Monodromy representations for generalized
	Knizhnik-Zamolodchikov equations and Hecke algebras},
	Publ. Res. Inst. Math. Sci. {\bf 27} (1991), 711--726.
\bibitem[Che4]{Ch4} I.~Cherednik,
	{\em Integration of quantum many-body problems by affine
	Knizhnik-Zamolodchikov equations},
	preprint RIMS {\bf 776} (1991),
	Adv. Math. {\bf 106} (1994), 65--95.
\bibitem[Che5]{Ch5} I.~Cherednik,
	{\em Double affine Hecke algebras, Knizhnik-Zamolodchikov equations,
	and Macdonald's operators},
	Internat. Math. Res. Notices {\bf 9} (1992), 171--180.
\bibitem[Che6]{Ch} I.~Cherednik,
	{\em Double affine Hecke algebras and difference Fourier transforms},
	 Invent. Math. {\bf 152} (2003), 213--303.
\bibitem[Che7]{Chl} I.~Cherednik,
	``Double affine Hecke algebras'',
	Cambridge University Press, 2005.
\bibitem[CheMa]{ChMa} I.~Cherednik and Y.~Markov,
	{\em Hankel transform via double Hecke algebra},
	in ``Iwahori-Hecke algebras and their representation theory'',
	1--25, Lecture Notes in Math., 1804, Springer Verlag, 2002. 
\bibitem[Chm]{Chm} T.~Chmutova,
	{\em Representations of the rational Cherednik algebras
	of dihedral type},
	preprint math.RT/0405383, to appear in Journal of Alg.
\bibitem[ChmEt]{ChmEt} T.~Chmutova and P.~Etingof,
	{\em On some representations of the rational Cherednik algebra},
	Represent. Theory {\bf 7} (2003), 641--650.
\bibitem[CrBoHo]{CrBoHo} W.~Crawley-Boevey and M.~Holland,
	{\em Noncommutative deformations of Kleinian singularities},
	Duke Math. J. {\bf 92} (1998), 605--635.
\bibitem[De]{De} C.~Dez\'el\'ee,
	{\em Repr\'esentations de dimension finie de l'alg\`ebre de
	Cherednik rationnelle},
	Bull. Soc. Math. France {\bf 131} (2003), 465--482.
\bibitem[Dr]{Dr} V.~Drinfeld,
	{\em Degenerate affine Hecke algebras and Yangians}
	Funktsional. Anal. i Prilozhen. {\bf 20} (1986), 69--70. 
\bibitem[Du1]{Du0} C.F.~Dunkl,
	{\em Differential-difference operators associated to reflection
	groups},
	Trans. Amer. Math. Soc. {\bf 311} (1989), 167--183.
\bibitem[Du2]{Duh} C.F.~Dunkl,
	{\em  Hankel transforms associated to finite reflection groups},
	in ``Hypergeometric functions on domains of positivity, Jack
	polynomials, and applications'' (Tampa, FL, 1991), 123--138,
	Contemp. Math., 138, Amer. Math. Soc., 1992. 
\bibitem[Du3]{Du1} C.F.~Dunkl,
	{\em Differential-difference operators and monodromy
	representations of Hecke algebras},
	Pacific J. Math. {\bf 159} (1993), 271--298. 
\bibitem[Du4]{Du} C.F.~Dunkl,
	{\em Singular polynomials and modules for the symmetric groups},
	preprint math.RT/0501494.
\bibitem[DuDeJOp]{DuDeJOp} C.F.~Dunkl, M.F.E.~De Jeu and E.~Opdam,
	{\em Singular polynomials for finite reflection groups},
	Transactions Amer. Math. Soc. {\bf 346} (1994), 237--256.
\bibitem[DuOp]{DuOp} C.F.~Dunkl and E.~Opdam,
        {\em Dunkl operators for complex reflection groups},
        Proc. London Math. Soc. {\bf 86} (2003), 70--108.
\bibitem[Et]{Et} P.~Etingof,
	{\em Cherednik and Hecke algebras of varieties with a finite group
	action},
	preprint math.QA/0406499(v3).
\bibitem[EtGanGi]{EtGaGi} P.~Etingof, W.~Liang Gan, and V.~Ginzburg,
	{\em Continuous Hecke algebras},
	preprint math.QA/0501192(v2).
\bibitem[EtGi]{EtGi} P.~Etingof and V.~Ginzburg,
        {\em Symplectic reflection algebras, Calogero-Moser space and
        deformed Harish-Chandra homomorphism}, Inv. Math. {\bf 147}
        (2002), 243--348.
\bibitem[EtMo]{EtMo} P.~Etingof and S.~Montarani,
	{\em Finite dimensional representations of symplectic reflection
	algebras associated to wreath products},
	preprint math.RT/0403250.
\bibitem[EtSt]{EtSt} P.~Etingof and E.~Strickland,
	{\em Lectures on quasi-invariants of Coxeter groups and
	the Cherednik algebra},
	Enseign. Math. {\bf 49} (2003), 35--65.
\bibitem[Gan]{Ga} W.L.~Gan,
	{\em Reflection functors and symplectic reflection algebras for
	wreath products},
	preprint math.RT/0502035.
\bibitem[GanGi]{GaGi} W.L.~Gan and V.~Ginzburg,
	{\em Almost-commuting variety, D-modules, and Cherednik Algebras},
	preprint math.RT/0409262(v2).
\bibitem[GarGr]{GaGr} H.~Garland and I.~Grojnowski,
	{\em Affine Hecke algebras associated to Kac-Moody groups},
	preprint q-alg/9508019.
\bibitem[Gi]{Gi} V.~Ginzburg,
	{\em On primitive ideals},
	Selecta Math. {\bf 9} (2003), 379--407.
\bibitem[GGOR]{GGOR} V.~Ginzburg, N.~Guay, E.~Opdam, and
	R.~Rouquier,
        {\em On the category $\mathcal{O}$ for rational Cherednik algebras},
        Inventiones Math. {\bf 154} (2003), 617--651
\bibitem[GiKal]{GiKa} V.~Ginzburg and D.~Kaledin,
	{\em Poisson deformations of symplectic quotient singularities},
	 Adv. Math. {\bf 186} (2004), 1--57.
\bibitem[GiKapVas]{GiKaVa} V.~Ginzburg, M.~Kapranov, and E.~Vasserot,
	{\em  Residue construction of Hecke algebras},
	Adv. Math. {\bf 128} (1997), 1--19.
\bibitem[Go1]{Go1} I.~Gordon,
	{\em Baby Verma modules for rational Cherednik algebras},
	Bull. London Math. Soc. {\bf 35} (2003), 321--336.
\bibitem[Go2]{Go2} I.~Gordon,
	{\em  On the quotient ring by diagonal invariants},
	Invent. Math. {\bf 153} (2003), 503--518.
\bibitem[GoSm]{GoSm} I.~Gordon and S.P.~Smith,
	{\em Representations of symplectic reflection algebras and
	resolutions of deformations of symplectic quotient singularities},
	 Math. Ann. {\bf 330} (2004), 185--200.
\bibitem[GoSt1]{GoSt1} I.~Gordon and J.T.~Stafford,
        {\em Rational Cherednik algebras and Hilbert schemes},
        preprint math.RA/0407516.
\bibitem[GoSt2]{GoSt2} I.~Gordon and J.T.~Stafford,
        {\em Rational Cherednik algebras and Hilbert schemes II:
        representations and sheaves}, preprint math.RT/0410293.
\bibitem[Gu1]{Gu} N.~Guay,
	{\em Projective modules in the category $\CO$ for the
	Cherednik algebra},
	J. Pure Appl. Algebra {\bf 182} (2003), 209--221.
\bibitem[Gu2]{Gu2} N.~Guay,
	{\em Cherednik algebras and Yangians},
	preprint, April 2005.
\bibitem[Kap]{Ka} M.~Kapranov,
	{\em  Double affine Hecke algebras and 2-dimensional local fields},
	J. Amer. Math. Soc. {\bf 14} (2001), 239--262.
\bibitem[La]{La} F.~Latour,
	{\em Representations of rational Cherednik algebras of rank one
	in positive characteristic},
	J. Pure Appl. Algebra {\bf 195} (2005), 97--112.
\bibitem[LeThi]{LeThi} B.~Leclerc and J.-Y.~Thibon,
	{\em Canonical bases of $q$-deformed Fock spaces},
	Internat. Math. Res. Notices {\bf 9} (1996), 447--456.
\bibitem[Ma]{Ma} M.~Martino,
	{\em The associated variety of a Poisson prime ideal},
	preprint math.RT/0405253.
\bibitem[Mo1]{Mo1} S.~Montarani,
	{\em On some finite dimensional representations of symplectic
	reflection algebras associated to wreath products},
	preprint math.RT/0411286(v2).
\bibitem[Mo2]{Mo2} S.~Montarani,
	{\em Finite dimensional representations of symplectic
	reflection algebras associated to wreath products II},
	preprint math.RT/0501156.
\bibitem[Mu]{Mu} I.M.~Musson,
	{\em Noncommutative Deformations of Type A Kleinian Singularities
	and Hilbert Schemes},
	preprint math.RT/0504543.
\bibitem[Op]{Op} E.~Opdam,
	{\em Bessel functions and the discriminant of a finite Coxeter group},
	Compositio Math. {\bf 85} (1993), 333--373. 
\bibitem[Rou]{Rou} R.~Rouquier,
	{\em $q$-Schur algebras and complex reflection groups, I},
	in preparation.
\bibitem[Su]{Su} T.~Suzuki,
	{\em Rational and trigonometric degeneration of the double
	affine Hecke algebra of type $A$},
	preprint math.RT/0502534.
\bibitem[VarVas1]{VarVas1} M.~Varagnolo and E.~Vasserot,
	{\em Schur duality in the toroidal setting},
	Comm. Math. Phys. {\bf 182} (1996), 469--483.
\bibitem[VarVas2]{VarVas} M.~Varagnolo and E.~Vasserot,
	{\em From double affine Hecke algebras to quantized affine Schur
	algebras},
	Int. Math. Res. Not. {\bf 26} (2004), 1299--1333.
\bibitem[Vas]{Va} E.~Vasserot,
	{\em Induced and simple modules of double affine Hecke algebras},
	 Duke Math. J. {\bf 126} (2005), 251--323.
\end{thebibliography}
\end{document}